\input amssym.def
\input amssym.tex

\magnification=1200

\catcode`\@=11

\hsize=120 mm   \vsize =175mm \hoffset=4mm    \voffset=10mm
\pretolerance=500 \tolerance=1000 \brokenpenalty=5000

%
%
%
%
%
%

\newif\ifpagetitre      \pagetitretrue
\newtoks\hautpagetitre  \hautpagetitre={\hfil}
\newtoks\baspagetitre   \baspagetitre={\hfil}

\newtoks\auteurcourant  \auteurcourant={\hfil}
\newtoks\titrecourant   \titrecourant={\hfil}

\newtoks\hautpagegauche \newtoks\hautpagedroite
\hautpagegauche={\hfil\the\auteurcourant\hfil}
\hautpagedroite={\hfil\the\titrecourant\hfil}

\newtoks\baspagegauche  \baspagegauche={\hfil\tenrm\folio\hfil}
\newtoks\baspagedroite  \baspagedroite={\hfil\tenrm\folio\hfil}

\headline={\ifpagetitre\the\hautpagetitre
\else\ifodd\pageno\the\hautpagedroite
\else\the\hautpagegauche\fi\fi}

\footline={\ifpagetitre\the\baspagetitre
\global\pagetitrefalse
\else\ifodd\pageno\the\baspagedroite
\else\the\baspagegauche\fi\fi}

\hautpagetitre={\hfill\tenbf preliminary version: Not for diffusion\hfill}
\hautpagetitre={\hfill}
\hautpagegauche={\tenrm\folio\hfill\the\auteurcourant}
\hautpagedroite={\tenrm\the\titrecourant\hfill\folio}
\baspagegauche={\hfil} \baspagedroite={\hfil}
\auteurcourant{Marmi, Moussa, Yoccoz}
\titrecourant{preliminary version: Not for diffusion}
\def\mois{\ifcase\month\or January\or February\or March\or April\or
May\or June\or July\or August\or September\or October\or November\or
December\fi}
\def\Date{\rightline{\mois\ /\ \the\day\ /\/ \the\year}}
\hfuzz=0.3pt
\font\tit=cmb10 scaled \magstep1

\def\IM{\mathop{\Im m}\nolimits}
\def\RE{\mathop{\Re e}\nolimits}
\def\H{\Bbb H}
\def\R{\Bbb R}
\def\T{\Bbb T}
\def\Z{\Bbb Z}

\def\C{\Bbb C}
\def\N{{\Bbb N}}

\def\hH+{\hat{\H}^{+}}
\def\hHZ+{\widehat{\H^{+}/\Z}}

\def\PGL2Z{\hbox{PGL}\, (2,\Z)}
\def\GL2Z{\hbox{GL}\, (2,\Z)}
\def\SL2Z{\hbox{SL}\, (2,\Z)}

\def\cA{{\cal A}}

\def\cC{{\cal C}}
\def\cD{{\cal D}}

\def\cI{{\cal I}}

\def\cN{{\cal N}}

\def\cU{{\cal U}}

\def\cZ{{\cal Z}}

\def\proof{\noindent{\it Proof.\ }}
\def\qed{\hfill$\square$\par\smallbreak}
\def\Proc#1#2\par{\medbreak \noindent {\bf #1\enspace }{\sl #2}%
\par\ifdim \lastskip <\medskipamount \removelastskip \penalty 55\medskip \fi}%
\def\Def#1#2\par{\medbreak \noindent {\bf #1\enspace }{#2}%
\par\ifdim \lastskip <\medskipamount \removelastskip \penalty 55\medskip
\fi}
\def\qed{\hfill$\square$\par\smallbreak}
\def\hfl#1#2{\smash{\mathop{\hbox to 12mm{\rightarrowfill}}
\limits^{\scriptstyle#1}_{\scriptstyle#2}}}

\def\build#1_#2^#3{\mathrel{\mathop{\kern 0pt#1}\limits_{#2}^{#3}}}



\rightline{March 25, 2004}
\medskip
\centerline{\tit The cohomological equation for Roth type interval exchange maps }
\bigskip
\centerline{S. Marmi\footnote{$^1$}{Scuola Normale Superiore, Piazza dei Cavalieri 7, 56126
Pisa, Italy},
P. Moussa\footnote{$^2$}{Service de Physique Th\'eorique, CEA/Saclay,
91191 Gif-Sur-Yvette, France}
and J.-C. Yoccoz\footnote{$^3$}{Coll\`ege de France, 3, Rue d'Ulm,
75005 Paris }}
\vskip 2. truecm

\bigskip
\centerline{CONTENTS}
\smallskip{
\item{0.} Introduction
\item{1.}  The continued fraction algorithm for interval exchange
maps
\item{1.1} Interval exchnge maps
\item{1.2} The continued fraction algorithm
\item{1.3} Roth type interval exchange maps
\item{2.} The cohomological equation
\item{2.1} The theorem of Gottschalk and Hedlund
\item{2.2} Special Birkhoff sums
\item{2.3} Estimates for functions of bounded variation
\item{2.4} Primitives of functions of bounded variation
\item{3.} Suspensions of  interval exchange maps
\item{3.1} Suspension data
\item{3.2} Construction of a Riemann surface
\item{3.3} Compactification of $M_\zeta^*$
\item{3.4} The cohomological equation for higher smoothness
\item{4.} Proof of full measure for Roth type
\item{4.1} The basic operation of the algorithm for suspensions
\item{4.2} The Teichm\"uller flow
\item{4.3} The absolutely continuous
invariant measure
\item{4.4} Integrability of $\log\Vert Z_{(1)}\Vert$
\item{4.5} Conditions (b) and (c) have full measure
\item{4.6} The main step
\item{4.7} Condition (a) has full measure
\item{4.8} Proof of the Proposition

\item{Appendix A} Roth--type conditions in a concrete family of
i.e.m.\
\item{Appendix B} A non--uniquely ergodic i.e.m.\ satsfying condition (a)
\item{References}

 }

\vfill\eject

\vskip .5 truecm \noindent {\bf 0. Introduction}

\vskip 1. truecm
Let $\alpha$ be an irrational number, $(q_n)_{n\in\N}$ be the
sequence of the denominators of its continued fraction expansion
and $(a_n)_{n\in\N}$ be the sequence of its partial quotients.
Roth type irrationals have several equivalent arithmetical
characterizations:
\item{$\bullet$} in terms of the rate of approximation by rational
numbers: for all $\varepsilon >0$ there exists a positive constant
$C_\varepsilon$ such that $|q\alpha -p|\ge C_\varepsilon
q^{1+\varepsilon}$ for all rationals $p/q$;
\item{$\bullet$} in terms of the growth rate of the denominators
of the continued fraction: $q_{n+1}=\hbox{O}\,
(q_n^{1+\varepsilon})$ for all $\varepsilon >0$;
\item{$\bullet$} in terms of the growth rate of the partial
quotients: $a_{n+1}=\hbox{O}\,
(q_n^{\varepsilon})$ for all $\varepsilon >0$.

In addition to these purely arithmetical characterizations an
equivalent definition arises naturally in the study of the
cohomological equation associated to the rotation $R_\alpha\, :\,
x\mapsto x+\alpha$ on the circle $\T=\R/\Z$: $\alpha$ is of Roth
type if and only if for all $r,s\in\R$ with $r>s+1\ge 1$ and for
all functions $\Phi$ of class $\cC^r$ on $\T$ with zero mean
$\int_\T\Phi dx=0$ there exists a unique function $\Psi$ of class
$\cC^s$ on $\T$ and with zero mean such that $\Psi-\Psi\circ
R_\alpha=\Phi$.

The class of Roth type irrationals enjoys several nice properties:
by the celebrated theorem of Roth all algebraic irrationals are of
Roth type. Moreover the set of Roth type numbers has full measure
and is invariant under the natural action of the modular group
$\hbox{SL}\, (2,\Z )$.

The goal of this paper is to characterize a class of interval
exchange maps (i.e.m.'s) with similar properties (especially for
the solutions of the associated cohomological equation and the
fact of being a full measure class).

\vskip .5 truecm \noindent {\bf 0.1 Interval exchange maps}

\vskip .5 truecm \noindent
Let $\cal A$ denote an alphabet with
$d\ge 2$ elements. Let $I$ be an interval and $(I_\alpha
)_{\alpha\in\cA}$ a partition of $I$ into $d$ subintervals. An
interval exchange map $T$ is an invertible map of $I$ which is a
translation on each $I_\alpha$. Thus $T$ is orientation--preserving and
preserves Lebesgue measure.

When $d=2$ then $T$ is just a rotation (modulo identification of
the endpoints of $I$). It can be thought as the first return map
of a linear flow on a two--dimensional torus on a transversal
circle. Analogously when $d\ge 3$ by singular suspension any
i.e.m.\ is related to the linear flow on a suitable translation
surface (see, e.g.\ [V1] for details, or section 3 below)
typically having genus higher than $2$. A well--known dictionary
between translation surfaces and Riemann surfaces relates i.e.m.'s
to the theory of measured foliations on surfaces (see, e.g.\ [FLP]
for a introduction to measured foliations). Finally i.e.m.'s are
related to the study of rational polygonal billiards (see [Ar],
[Ta] and [KH], Chapter 14, for a general introduction to i.e.m.'s,
flows on surfaces and polygonal billiards).

Typical i.e.m.'s are minimal (this is guaranteed by a condition
due to Keane [Ke1] which is automatically dealt with if the
intervals lengths are rationally independent) but note that
ergodic properties of minimal i.e.m.'s can differ substantially
from those of circle rotations: first they need not be uniquely
ergodic [Ke2, KN, Co], and second, being ergodic they can be
weakly mixing [KS, V3,V4]. On the other hand uniquely ergodic
i.e.m.'s are generic [KR] and Keane's conjecture that almost every
i.e.m.\ is uniquely ergodic was proven independently by Masur and
Veech [Ma, V2], see also [Ker, Re].

One of the most important consequences for us of Keane's condition
is that it allows to introduce and to iterate indefinitely
continued fraction algorithms that generalize the classical
algorithm (corresponding to the choice $d=2$) [Ra, V2, Z1]. Both
the Rauzy--Veech continued fraction algorithm and its accelerated
version due to Zorich are ergodic w.r.t.\ an absolutely continous
invariant measure in the space of i.e.m.'s. However in the case of
the Rauzy--Veech continued fraction the measure has infinite mass
whereas the invariant measure for the Zorich algorithm has finite
mass. The ergodic properties of the continued fraction map and of
the related Teichm\"uller flow (see Section 4.2 for its
definition) have been studied in detail [V5, V6, V7, Z2, Z4, Fo2].

\vskip .5 truecm \noindent {\bf 0.2 The cohomological equation}

\vskip .5 truecm \noindent Our study of the cohomological equation
for i.e.m.'s has been prompted by Forni's [Fo1] celebrated paper
on the cohomological equation associated to linear flows on
surfaces of higher genus. Let us first state our main theorem.

We will denote $\hbox{BV}\, (\sqcup I_\alpha)$ (resp.\
$\hbox{BV}_*\, (\sqcup I_\alpha)$) the space of functions
$\varphi$ whose restriction to each of the intervals $I_\alpha$ is
a function of bounded variation (resp.\ the hyperplane of
$\hbox{BV}\, (\sqcup I_\alpha)$ made of functions whose integral
on the disjoint union $\sqcup I_\alpha$ vanishes). We will also
denote $\hbox{BV}^1_*\, (\sqcup I_\alpha)$) the space of functions
$\varphi$ which are absolutely continuous on each $I_\alpha$ and
whose first derivative belongs to $\hbox{BV}_*\, (\sqcup
I_\alpha)$.

Our first main result can be stated as follows:

\vskip .3 truecm\noindent {\bf Theorem A.}{\it Let $T$ be an interval
exchange map with
the Keane property and of Roth type. Let $\Phi\in\hbox{BV}^1_*\, (\sqcup
I_\alpha)$.  There exists a function $\chi$ constant on each
interval $I_\alpha$
and a bounded function $\Psi$ such that
}
$$\Psi
-\Psi\circ T=\Phi -\chi \; . $$

\vskip .5 truecm\noindent
To make the above statement precise we need to define Roth type i.e.m.'s. This is
the subject of section 1.3 below. For the time being we will
content ourselves with briefly describing the three conditions
which a Roth type i.e.m.\ must satisfy:
\item{(a)} The first condition is a growth rate condition for the
matrices appearing in an accelerated version of the Zorich
continued fraction algorithm (see Section 1.2.4 for details). This
condition is the precise analogue of the third of the equivalent
arithmetical characterizations of Roth type irrational numbers
given above.
\item{(b)} The second condition is a spectral condition which
guarantees unique ergodicity of Roth type i.e.m.'s. This condition
does not follow from condition (a) (see Appendix B for a
counterexample, and also [Ch]) but is automatically satisfied if
the i.e.m.\ is of {\it constant type} (i.e.\ the  matrices
considered in (a) have bounded norm).
\item{(c)} The third and last condition is a coherence condition.

The second main result of this paper is

\vskip .3 truecm\noindent {\bf Theorem B.}{\it Roth type interval
exchange maps form a full measure set in the space of all interval
exchange maps.}

\vskip .5 truecm\noindent Obviously, Theorem A is closely
connected to Forni's fundamental theorem [Fo1] on the
cohomological equation for area--preserving vector fields on
surfaces. By singular suspension (``zippered rectangles'', see
Section 3), one obtains from an interval exchange map an
area--preserving flow on a singular flat surface. Forni develops
some Fourier analysis tools in this context, which allows him to
solve the cohomological equation for almost every direction; our
methods are completely  different. He works in the Sobolev scale
and his methods allow to lose no more than $3+\varepsilon $
derivatives (for every $\varepsilon  >0$) [Fo3]. Our loss is
smaller and we get an explicit Diophantine condition. On the other
hand, given a singular flat surface, we do not know if almost
every direction leads to a Roth type interval exchange map.

The connection with singular flat surfaces explains the type of
regularity we introduce when we consider the cohomological
equation for more regular data: we still allow discontinuities for
$\Phi$ at the endpoints for the $I_\alpha$; on the other hand, we
require the solution $\Psi$ to be continuous on all of $I$. New
linear conditions on $\Phi$ appear by integration of the
cohomological equation. See Section 3 below for the precise
statements.

When the singular suspension of an i.e.m.\ $T$ is an invariant
foliation for a pseudo--Anosov diffeomorphism, the continued
fraction expansion of $T$ is eventually periodic. This implies a
strong version of condition (a).

Conditions (b) and (c) are also satisfied. Hence $T$ is of Roth
type (even of ``bounded type'') and Theorem A applies. This
answers positively a question raised by Forni ([Fo1], p.\ 342).

\vskip .5 truecm \noindent {\bf 0.3 Summary of the contents}

\vskip .5 truecm \noindent In the first section we introduce
interval exchange maps and we develop the continued fraction
algorithms to an extent which allows us to introduce Roth type
i.e.m.'s. The Keane property (see 1.1.6) does not only guarantee
that an i.e.m. is minimal but it also implies that the
Rauzy--Veech continued fraction algorithm (described in
1.2.1-1.2.3) can be iterated indefinitely. Accelerating the
Rauzy--Veech map by grouping together arrows with the same name in
the Rauzy diagram leads to the Zorich continued fraction algorithm
(described in 1.2.4) which has the advantage of having a finite
mass a.c.i.m.. On the other hand, since every name is taken
infinitely many times in the sequence of arrows in the Rauzy
diagram corresponding to a given i.e.m. one can produce a further
acceleration of the scheme by grouping together all arrows which
take all possible names but one: this leads to the algorithm we
will use in the definition of Roth type i.e.m.'s given in section
1.3 and already briefly described above. The notations and the
presentation of the Rauzy--Veech--Zorich algorithms follow closely
the expository paper [Y].

Section 2 is devoted to the study of the cohomological equation
and to the proof of our main theorem A.
When $T$ is a minimal homeomorphism of a compact space $X$,
we know from a theorem of
Gottschalk and Hedlund [GH] that a continuous function on $X$ is a
$T$--coboundary of some
continuous function as soon as its
Birkhoff sums at some point of $X$ are bounded (see Section 2.1.1).
An i.e.m.\ with
the Keane property is minimal but
not continuous. Nevertheless, a Denjoy-like construction (see Section 2.1.2)
allows to apply Gottschalk--Hedlund's theorem and
conclude that a continuous function whose Birkhoff sums at some point are bounded is the
$T$ -coboundary of a bounded function. The next step in the proof
is the reduction of the control of a general Birkhoff sum to the
control of those {\it special} Birkhoff sums which are obtained by
considering the return times of the point under iteration of the
map (Section 2.2). These can be conveniently analyzed using the
continued fraction. The estimates of these special Birkhoff sums
for functions of bounded
variation are given in Section 2.3 and the proof of the theorem is
completed in Section 2.4.

In Section 3 we first recall how to construct a linear flow on a
translation surface starting from an i.e.m.\ and certain
suspension data (Sections 3.1--3.3). Then we relate the discrete
cohomological equation for i.e.m.'s to the continuous one for the
vertical (area--preserving) vector field constructed by
suspension: this allows us to consider more regular data
(i.e.\ belonging to the space
$\hbox{BV}^r_*$ of functions whose $r$--th derivative has
bounded variation on each $I_\alpha$ and all intermediate derivatives have
zero mean on $\sqcup I_\alpha$ ).
We prove that for those the loss of differentiability in solving the
cohomological equation is the same as for functions in
$\hbox{BV}^1_*$ (Section 3.4).

Section 4 is devoted to the proof of theorem B, i.e.\ that Roth
type i.e.m.'s have full measure. To this purpose we need to
describe how the Rauzy--Veech map acts at the level of the
suspension data (Section 4.1). Then we combine the continued
fraction algorithm (in Zorich form) with the Teichm\"uller flow in
order to get a version which is normalized w.r.t.\ scales (Section
4.2). A careful comparison between the a.c.i.m.\ for the continued
fraction map and the Lebesgue measure is carried out in Section
4.3 whereas in Section 4.4, following Zorich [Z1] we prove the
integrability condition on the matrices needed to apply Oseledets
multiplicative ergodic theorem. Then conditions (b) and (c) in the
definition of Roth type i.e.m.'s have full measure (Section 4.5)
by  Oseledets theorem and the almost sure existence of a spectral
gap proved by Veech in [V3]. Showing that condition (a) also has
full measure requires more work and more precise informations on
the combinatorics of the continued fraction map. This is
summarized in a Proposition stated in Section 4.6 and proved in
Section 4.8 whereas in Section 4.3 we show how to conclude the
proof of theorem B by putting together the results of Sections 4.3
and 4.7 and applying a Borel--Cantelli argument.

The two appendices are devoted to the construction of concrete
examples of Roth type i.e.m.'s and to the construction of
non--uniquely ergodic i.e.m.'s satisfying condition (a) in Roth
type (but of course not condition (b)).

\vskip 1. truecm \noindent {\bf Acknowledgements} We are grateful
to G. Forni for many stimulating discussions.  This research has
been supported by the  following institutions: CNR, CNRS, MURST,
INDAM, the French--Italian University, the Coll\`ege de France and
the Scuola Normale Superiore. We are also grateful to the two
former institutions and to the Centro di Ricerca Matematica
``Ennio De Giorgi'' in Pisa for hospitality.

\vfill\eject \noindent {\bf 1. The continued fraction algorithm
for interval exchange maps}

\vskip .5 truecm \noindent {\bf 1.1 Interval exchange maps}

\vskip .5 truecm \noindent {\bf 1.1.1}

\vskip .3 truecm\noindent  An interval exchange map (i.e.m.\ ) is
determined by combinatorial data on one side, length data on the
other side.

The combinatorial data consists of a finite set $\cal A$ of names
for the intervals and of two bijections $(\pi_0, \pi_1)$ from
$\cA$ onto $\{1, \ldots ,d\}$ (where $d$ is the cardinality of
$\cA$): these indicate in which order the intervals are met before
and after the map.

The length data  $(\lambda_\alpha)_{\alpha\in\cA}$ give the length
$\lambda_\alpha >0$ of the corresponding interval. More precisely,
we set $$\eqalign{ I_\alpha & := [0,\lambda_\alpha
)\times\{\alpha\}\, , \cr \lambda^* & :=
\sum_{\alpha\in\cA}\lambda_\alpha\, , \cr I & := [0,\lambda^* )\,
.\cr} $$

We then define, for $\varepsilon =0,1$, a bijection
$j_\varepsilon$ from $\sqcup_{\alpha\in\cA}I_\alpha$ onto $I$: $$
j_\varepsilon (x,\alpha ) =  \sum_{\pi_\varepsilon(\beta
)<\pi_\varepsilon (\alpha )} \lambda_\beta\; . $$

The i.e.m.\ $T$ associated to these data is the bijection
$T=j_1\circ j_0^{-1}$ of $I$.

\vskip .5 truecm \noindent {\bf 1.1.2} If $\cA , \pi_0 , \pi_1 ,
\lambda_\alpha $ are as above and $\chi\, : \, \cA'\rightarrow\cA$
is a bijection, we can define a new set of data by $$\eqalign{
\pi_\varepsilon ' & = \pi_\varepsilon\circ\chi\;\; , \;
\varepsilon = 0,1\; , \cr \lambda_{\alpha '}' & = \lambda_{\chi
(\alpha ')}\;\; , \; \alpha '\in\cA'\; .\cr} $$

Obviously, the new i.e.m.\ $T'$ determined by these data is the
same, except for names, than the old one. In particular, we could
restrict to consider {\it normalized combinatorial data}
characterized by $$\cA=\{1,\ldots ,d\}\; , \; \; \pi_0 =
\hbox{id}_\cA\; .$$

However, this leads to later to more complicated formulas in the
continued fraction algorithm because the basic operations on
i.e.m.\ 's do not preserve normalization.

\vskip .5 truecm \noindent {\bf 1.1.3} Given  combinatorial data
$(\cA , \pi_0 , \pi_1 )$, we set, for $\alpha , \beta \in \cA$ $$
\Omega_{\alpha ,\beta} = \cases{ +1 & if $\pi_0(\beta )>\pi_0
(\alpha )\; , \; \pi_1(\beta )<\pi_1 (\alpha )$, \cr -1 & if
$\pi_0(\beta )<\pi_0 (\alpha )\; , \; \pi_1(\beta )>\pi_1 (\alpha
)$, \cr 0 & otherwise. } $$ The matrix $\Omega =
\left(\Omega_{\alpha ,\beta}\right)_{(\alpha ,\beta )\in\cA^2}$ is
antisymmetric.

Let $(\lambda_\alpha)_{\alpha\in\cA}$ be the length data and let
$T$ be the associated i.e.m.\ . For $\alpha\in\cA$, $y\in
j_0(I_\alpha )$, we have $$ T(y)=y+\delta_\alpha\; , $$ where the
{\it translation vector} $\delta=( \delta_\alpha )_{\alpha\in\cA}$
is related to the {\it length vector} $\lambda = (\lambda_\alpha
)_{\alpha\in\cA}$ by: $$\delta = \Omega\lambda\; . $$

\vskip .5 truecm \noindent {\bf 1.1.4} There is a {\it canonical
involution} ${\cal I} $ acting on the set of combinatorial data
which exchange $\pi_0$ and $\pi_1$. For any set $(\lambda_\alpha
)_{\alpha\in\cA}$ of length data, the intervals $I_\alpha , I$ are
unchanged, but $j_0$ and $j_1$ are exchanged and $T$ is replaced
by  $T^{-1}$. The matrix  $\Omega$ is replaced by $-\Omega$ and
the translation vector  $\delta$ by $-\delta$. Observe that ${\cal
I} $ does not respect the combinatorial normalization.

\vskip .5 truecm \noindent {\bf 1.1.5} In the following, we will
always consider only combinatorial data $(\cA , \pi_0 , \pi_1 )$
which are {\it admissible}, meaning that for all $k=1,2,\ldots
,d-1$, we have $$ \pi_0^{-1}(\{1,\ldots ,k\})\not=
\pi_1^{-1}(\{1,\ldots ,k\})\; . $$

Indeed, if we had $\pi_0^{-1}(\{1,\ldots ,k\})=
\pi_1^{-1}(\{1,\ldots ,k\})$ for some $k<d$, for any length data
$(\lambda_\alpha )_{\alpha\in\cA}$, the interval $I$ would
decompose into two disjoint invariant subintervals and the study
of the dynamics would be reduced to simpler combinatorial data.

\vskip .5 truecm \noindent {\bf 1.1.6 The Keane property} Let $T$
be an i.e.m.\ defined by combinatorial data $(\cA , \pi_0 , \pi_1
)$ and length data $(\lambda_\alpha )_{\alpha\in\cA}$.

\vskip .3 truecm \noindent {\bf Definition} {\it A } connexion
{\it for $T$ is a triple $(\alpha ,\beta ,m)$ where $\alpha
,\beta\in\cA$, $\pi_0 (\beta )>1$, $m$ is a positive integer, and}
$$T^m (j_0(0,\alpha ))=j_0(0,\beta )\; . $$ {\it We say that $T$
has the} Keane property {\it if there is no connexion for $T$.}

\vskip .3 truecm\noindent It turns out that this property is the
appropriate notion of irrationality for i.e.m.\ . The following
results are due to Keane ([Ke1]):
\item{$\bullet$} An i.e.m.\ with Keane's property is minimal
(i.e.\ all orbits are dense).
\item{$\bullet$} If the length data are rationally independent
(and the combinatorial data are admissible) then $T$ has Keane's
property.

\vskip 1. truecm \noindent {\bf 1.2 The continued fraction
algorithm}

\vskip .5 truecm \noindent {\bf 1.2.1 The basic operation. (Rauzy
[Ra], Veech [V2]) } Let $T$ be an i.e.m.\ defined by combinatorial
data $(\cA , \pi_0 , \pi_1 )$ and length data $(\lambda_\alpha
)_{\alpha\in\cA}$. We assume as always that the combinatorial data
are admissible.

We denote by $\alpha_0, \alpha_1$ the (distinct) elements of
$\cA$ such that $$\pi_0(\alpha_0)=\pi_1(\alpha_1)=d\; .$$ Observe
that if $\lambda_{\alpha_0}=\lambda_{\alpha_1}$, the triple
$(\alpha_0,\alpha_1,1)$ is a connexion and $T$ has not the Keane
property.

We now assume that $\lambda_{\alpha_0}\not=\lambda_{\alpha_1}$ and
define $\varepsilon\in\{0,1\}$ by $$ \lambda_{\alpha_\varepsilon}
= \hbox{Max}\, (\lambda_{\alpha_0} ,\lambda_{\alpha_1})\; . $$ We
set $$\eqalign{ \hat{\lambda}^* & =
\lambda^*-\lambda_{\alpha_{1-\varepsilon}}\; , \cr \hat{I} & = [0,
\hat{\lambda}^* )\subset I\; , \cr} $$ and define $\hat{T}\, :\,
\hat{I}\rightarrow \hat{I}$ to be the first return map of $T$ in
$\hat{I}$.

When $\varepsilon =0$ we have $$ \hat{T}(y) = \cases{ T(y) \; & if
$ y\notin j_0(I_{\alpha_1})$, \cr T^2(y) \; & if $y\in
j_0(I_{\alpha_1})$. \cr}$$ When $\varepsilon =1$ we have similarly
$$ \hat{T}^{-1} (y) = \cases{ T^{-1}(y) \; & if $ y\notin
j_1(I_{\alpha_0})$, \cr T^{-2}(y) \; & if $y\in
j_1(I_{\alpha_0})$. \cr}$$

In both cases, it appears that $\hat{T}$ is again an interval
exchange map which can be defined using the same alphabet $\cA$.
The length data for $\hat{T}$ are given by $$ \eqalign{
\hat{\lambda}_\alpha & = \lambda_\alpha \;\;\hbox{if}\;
\alpha\not=\alpha_\varepsilon\; ,\cr
\hat{\lambda}_{\alpha_\varepsilon} &= \lambda_{\alpha_\varepsilon}
- \lambda_{\alpha_{1-\varepsilon}} \; . \cr} $$

The combinatorial data $(\hat{\pi}_0,\hat{\pi}_1)$ for $\hat{T}$
are given by $$\hat{\pi}_\varepsilon =\pi_\varepsilon$$ and
$$\hat{\pi}_{1-\varepsilon} (\alpha) = \cases{
\pi_{1-\varepsilon}(\alpha ) & if $\pi_{1-\varepsilon}(\alpha )\le
\pi_{1-\varepsilon}(\alpha_\varepsilon )$,\cr
\pi_{1-\varepsilon}(\alpha ) +1 & if
$\pi_{1-\varepsilon}(\alpha_\varepsilon )
<\pi_{1-\varepsilon}(\alpha )<d$,\cr
\pi_{1-\varepsilon}(\alpha_\varepsilon ) +1 & if
$\pi_{1-\varepsilon}(\alpha )=d$.\cr }$$

We rewrite the relation between old and new length data as $$
\lambda = V\hat{\lambda}\; , $$ where $$V=\hbox{\bf
I}+E_{\alpha_\varepsilon\alpha_{1-\varepsilon}}$$ has now non
negative integer coefficients and belongs to the group
$\hbox{SL}\,(\Z^\cA)$. We also write $$
(\hat{\pi}_0,\hat{\pi}_1)=R_\varepsilon (\pi_0,\pi_1)$$ and observe
that these new combinatorial data are admissible.

\vskip .5 truecm \noindent {\bf 1.2.2 Rauzy diagrams } Let $\cA$
be an alphabet. We define an oriented graph as follows. The
vertices are the admissible pairs $(\pi_0, \pi_1)$. Each vertex
$(\pi_0,\pi_1)$ is the starting point of exactly two arrows with
endpoints at $R_0 (\pi_0, \pi_1)$ and $R_1 (\pi_0, \pi_1)$. The
arrow connecting $(\pi_0,\pi_1)$ to $R_\varepsilon (\pi_0,\pi_1)$ is
said to be of {\it type} $\varepsilon$.

The operations $R_0,R_1$ are obviously invertible. Therefore each
vertex is also the endpoint of exactly two arrows, one of each
type.

To each arrow in the graph, we associate a {\it name} in $\cA$: it
is the element $\alpha_\varepsilon$ such that $\pi_\varepsilon
(\alpha_\varepsilon)=d$ (where $(\pi_0,\pi_1)$ is the starting
point of the arrow and $\varepsilon$ is its type). The element
$\alpha_{1-\varepsilon}$ will then be called the {\it secondary
name} of this arrow.

A {\it Rauzy diagram} is a connected component of this oriented
graph.

Obviously, the {\it Rauzy operations} $R_0,R_1$ commute with
change of names (see 1.2). Up to change of names, there is only
one Rauzy diagram with $d=\hbox{card}\,\cA=2$, and one with
$d=\hbox{card}\,\cA=3$.


In the diagrams
in figure 1 the pair $(\pi_0,\pi_1)$ is denoted by the symbol
$\matrix{\pi_0^{-1}(1) & \ldots & \pi_0^{-1}(d)\cr \pi_1^{-1}(1) &
\ldots & \pi_1^{-1}(d)\cr}$. For $d=\hbox{card}\,\cA=4$ there are
$2$ distinct Rauzy diagrams: (see figure 2).


In each of these diagrams, the symmetry with respect to the
vertical axis corresponds to the action of the canonical
involution.

In the last diagram, there is a further symmetry with respect to
the center of the diagram, which corresponds to the exchange of
the names $B_0,B_1$. This is a monodromy phenomenon: to each
admissible pair $(\pi_0,\pi_1)$, one can associate the permutation
$\pi :=\pi_1\circ\pi_0^{-1}$ of $\{1,\ldots ,d\}$, which is
invariant under change of names. When we identify vertices with
the same permutation, we obtain a {\it reduced Rauzy diagram} and
we have a covering map from the Rauzy diagram onto the reduced
Rauzy diagram.

In the first three examples above, the covering map is an
isomorphism. In the last example, the degree of the covering map
is $2$ and the reduced Rauzy diagram is given in figure 3
%
%
where $\pi$ is denoted by $(\pi^{-1}(1),\ldots ,\pi^{-1}(d))$.

\vskip .5 truecm \noindent {\bf 1.2.3 The Rauzy--Veech algorithm}
Let $T$ be an i.e.m.\ with admissible combinatorial data. If $T$
has Keane's property, the basic operation is defined for $T$ and
it is immediate to check that the new i.e.m.\ $\hat{T}$ again has
Keane's property. Therefore we can iterate the basic operation and
generate a sequence $(T^{(n)})_{n\ge 0}$ of i.e.m.\ 's (with
$T^{(0)}=T$). We will denote $(\pi_0^{(n)}, \pi_1^{(n)})$ the
combinatorial data of $T^{(n)}$, by
$(\lambda_\alpha^{(n)})_{\alpha\in\cA}$ its length data, by
$\gamma^{(n)}$ the arrow in the Rauzy diagram connecting
$(\pi_0^{(n-1)},\pi_1^{(n-1)})$ to $(\pi_0^{(n)},\pi_1^{(n)})$, by
$V^{(n)}$ the matrix relating $\lambda^{(n-1)}$ to $\lambda^{(n)}$
through $$ \lambda^{(n-1)}=V^{(n)}\lambda^{(n)}\; . $$

Conversely, it is not difficult to check that when $T$ has a
connexion, the algorithm has to stop because one runs at some
point in the equality case $\lambda_{\alpha_0}=\lambda_{\alpha_1}$
in the basic operation.

\vskip .3 truecm \noindent {\bf Proposition}{\it Each name in
$\cA$ is taken infinitely many times by the sequence of arrows
$(\gamma^{(n)})_{n>0}$.}

\vskip .3 truecm\noindent\proof Let $\cA'$ be the set of names
which are taken infinitely many times and let
$\cA"=\cA\setminus\cA'$. Replacing $T$ by some $T^{(N)}$, we can
assume that names in $\cA"$ are not taken at all. Then the lengths
$\lambda_\alpha^{(n)}$, $\alpha\in\cA"$, do not depend on $n$. But
then elements $\alpha\in\cA"$ can only appear as secondary names
at most finitely many times. Replacing again $T$ by some
$T^{(N)}$, we can assume that secondary names are never in $\cA"$.
Then the sequences $(\pi_\varepsilon^{(n)}(\alpha))_{n>0}$, for
$\varepsilon\in\{0,1\}$, $\alpha\in\cA"$, are non decreasing and
we can assume (replacing once again $T$ by some $T^{(N)}$) that
they are constant.

We now claim that we must have $\pi^{(0)}_\varepsilon(\alpha ")<
\pi^{(0)}_\varepsilon(\alpha ')$ for all $\alpha "\in\cA"$,
$\alpha '\in\cA'$ and $\varepsilon\in\{0,1\}$. Because the pair
$(\pi_0^{(0)},\pi_1^{(0)})$ is admissible, this implies
$\cA'=\cA$. To prove the claim, assume that there exist $\alpha
'\in\cA'$, $\alpha "\in\cA"$, $\varepsilon\in\{0,1\}$ with
$\pi^{(0)}_\varepsilon(\alpha ')< \pi^{(0)}_\varepsilon(\alpha
")$.  As $\pi^{(n)}_\varepsilon(\alpha ")=
\pi^{(0)}_\varepsilon(\alpha ")$ for all $n\ge 0$, we can never
have $\pi^{(n)}_\varepsilon(\alpha ')=d$ for some $n>0$. By
definition of $\cA'$, there must exist $n\ge 0$ such that
$\pi^{(n)}_{1-\varepsilon}(\alpha ')=d$; but then
$\pi^{(n+1)}_\varepsilon(\alpha ")\not=
\pi^{(0)}_\varepsilon(\alpha ")$, which gives a contradiction.
\qed

\vskip .5 truecm\noindent {\bf 1.2.4 The Zorich algorithm and its
accelerations}

When $d=2$, setting $x=\lambda_B/\lambda_A$, the basic operation
reduces to the well--known map $$ g(x)=\cases{ {x\over 1-x}\; &
for $0<x<1/2$, \cr {1-x\over x}\; & for $1/2<x<1$,\cr}$$ with a
parabolic fixed point at $0$. There is a unique absolutely
continuous invariant measure, namely $dx/x$, but this measure is
infinite. On the other hand, the Gauss map generating the
continued fraction algorithm has $dx/(1+x)$ as a finite a.c.i.m.\
.

For i.e.m.'s with more intervals, identifying i.e.m.'s with
proportional length data (and the same combinatorial data), Veech
has shown [V2] that there exists again for the basic operation a
unique absolutely continuous invariant measure. Again this measure
is infinite. Zorich has discovered ([Z1]) how to concatenate
several steps of the basic operations in order to get a finite
a.c.i.m.\ .

Let $T$ be an i.e.m.\ with Keane's property, $T^{(n)}$,
$\gamma^{(n)}$, $V^{(n)}$ the data generated by the iteration of
the basic operation. Let also $1\le D<d$. We define inductively an
increasing sequence $n_D(k)=n_D(k,T)$ by setting $n_D(0)=0$ and:

\noindent $n_D(k+1)$ is the largest integer such that no more than
$D$ names are taken by the $\gamma^{(n)}$, for $n_D(k)<n\le
n_D(k+1)$.

The sequence is well defined because of the Proposition above.

Obviously, for $1<D<d$, $(n_D(k))_{k\ge 0}$ is a subsequence of
$(n_{D-1}(l))_{l\ge 0}$.

We will define, for $k>0$ $$ Z_{(D)}(k)=V^{(n_D(k-1)+1)}\cdots
V^{(n_D(k))}\; . $$

The case $D=1$ is the one considered by Zorich ([Z1]). We will on
the other hand be interested in the case $D=d-1$.

When the context is clear, we will simply write $Z(k)$ for
$Z_{(d-1)}(k)$ and $T^{(k)}$ for $T^{(n_{d-1}(k))}$,
$\lambda^{(k)}$ for its length data. With these notations, we have
$$ \lambda^{(k)} = Z(k+1)\lambda^{(k+1)}\; . $$ We will also set,
for $k<l$ $$ Q(k,l)=Z(k+1)\cdots Z(l)\; , $$ in order to have
$$\lambda^{(k)}=Q(k,l)\lambda^{(l)}\; . $$ We will also write
$Q(l)$ for $Q(0,l)$. The coefficients $Q_{\alpha\beta}(k,l)$ have
the following interpretation. Let $I^{(k)}=\sqcup_{\alpha\in\cA}
j_0(I_\alpha^{(k)})$ be the domain of $T^{(k)}$. For $l\ge k$, we
have $I^{(l)}\subset I^{(k)}$ and $T^{(l)}$ is the first return
map of $T^{(k)}$ in $I^{(l)}$. Then, the non negative integer
$Q_{\alpha\beta}(k,l)$ is the time spent in $j_0(I_\alpha^{(k)})$
by any point of $j_0(I_\beta^{(l)})$ until it returns in
$I^{(l)}$.

We will also introduce $$Q_\beta
(k,l)=\sum_{\alpha\in\cA}Q_{\alpha\beta}(k,l)\; , $$ which is the
return time in $I^{(l)}$ for points in $I_\beta^{(l)}$.

The following Lemma is the main reason to choose $D=d-1$ rather
than $D=1$.

\vskip .3 truecm\noindent {\bf Lemma.}{\it Let $T$ satisfy Keane's
condition. Assume that $$l\ge \cases{ k+2d-3 &if $d\ge 3$\cr k+2
&if $d=2$.\cr}$$ Then, for all $\alpha ,\beta\in\cA$, we have
$Q_{\alpha\beta}(k,l)>0$.}

\vskip .3 truecm\noindent \proof Replacing $T$ by $T^{(k)}$, it is
sufficient to consider the case $k=0$. For $r\ge 0$, set $$
\hat{Q}(r)=V^{(1)}\cdots V^{(r)}\; ; $$ as the diagonal terms of
the $V$ matrices are equal to $1$ (and all the terms are non
negative) we have $$\hat{Q}_{\alpha\beta}(r)>0\Rightarrow
\hat{Q}_{\alpha\beta}(r+1)>0\; . $$

Fix $\alpha , \beta\in\cA$. We will construct a sequence of {\it
distinct} indices $\alpha_1=\alpha,\alpha_2,\ldots
,\alpha_s=\beta$ and integers $r_1=0<r_2<\ldots <r_s$ such that
$$\hat{Q}_{\alpha_1\alpha_j}(r)>0\;\;\hbox{for}\;r\ge r_j\; . $$
If $\alpha=\beta$, $s=1$, $r_1=0$ and the property is satisfied.
Otherwise, let $r_2$ be the smallest positive integer such that
the name of $\gamma^{(r_2)}$ is $\alpha_1$, and let $\alpha_2$ be
the secondary name of $\gamma^{(r_2)}$; we have
$\alpha_2\not=\alpha_1$ and $V^{(r_2)}_{\alpha_1\alpha_2}=1$ hence
$\hat{Q}_{\alpha_1\alpha_2}(r)>0$ for $r\ge r_2$.

Assume that $\alpha_1,\ldots ,\alpha_j,r_1,\ldots ,r_j$ have been
constructed, with $\beta\not=\alpha_l$ for $1\le l\le j$. Let
$r_j'$ be the smallest integer $>r_j$ such that the name of
$\gamma^{(r_j')}$ does not belong to $\{\alpha_1,\ldots
,\alpha_j\}$ and let $r_{j+1}$ be the smallest integer $>r_j'$
such that the name of $\gamma^{(r_{j+1})}$ belongs to
$\{\alpha_1,\ldots ,\alpha_j\}$; let $\alpha_{j+1}$ be the
secondary name of $\gamma^{(r_{j+1})}$. Then $\alpha_{j+1}$ is the
name of $\gamma^{(r_{j+1}-1)}$ and therefore is distinct from
$\alpha_1,\ldots ,\alpha_j$. By construction, we have, for some
$1\le l\le j$ $$V^{(r_{j+1})}_{\alpha_l\alpha_{j+1}}=1\; , $$ and
also $$\hat{Q}_{\alpha_1\alpha_l}(r_{j+1}-1)>0$$ because
$r_{j+1}>r_l$. We conclude that
$$\hat{Q}_{\alpha_1\alpha_{j+1}}(r)>0\;\;\hbox{for}\; r\ge
r_{j+1}\; . $$ At some point we will obtain $\alpha_s=\beta$. It
remains to see how many steps of the accelerated Zorich algorithm
(with $D=d-1$) are needed to attain $r_s$. Obviously, we have
$r_2\le n_{d-1}(1)+1$. Then, for $2\le j<d-1$, we have $$\eqalign{
r_j' & \le n_{d-1}(2j-2)\; ,\cr r_{j+1} &\le  n_{d-1}(2j-1)\;
.\cr}$$ Finally, when $s=d>2$, we have $$\eqalign{r_{d-1}' &
\le n_{d-1}(2d-4)+1\; ,\cr r_{d} &\le  n_{d-1}(2d-3)\; .\cr}$$
\qed

\vskip 1. truecm \noindent {\bf 1.3 Roth--type interval exchange
maps}

\vskip .5 truecm\noindent Roth--type i.e.m.\ should satisfy
Keane's condition so that the continued fraction algorithm is
defined, and three further conditions which are now explained.

\vskip .5 truecm \noindent {\bf 1.3.1 Size of the $Z$ matrices}

\vskip .3 truecm\noindent Take $D=d-1$ in 1.2.4. We will first ask
for the $Z$ matrices to be not too big in the following sense:
\item{(a)} {\it for every $\varepsilon >0$ there exists $C_\varepsilon
>0$ such that for all $k\ge 0$ we have}
$$ \Vert Z(k+1)\Vert \le C_\varepsilon\Vert
Q(k)\Vert^\varepsilon\; . $$

When $d=2$, this amounts exactly to the classical Roth type
approximation property for an irrational number $\theta$: for all
$\varepsilon >0$, there exists $\gamma_\varepsilon >0$ such that
for all rational $p/q$ one has $$|\theta-p/q|\ge\gamma_\varepsilon
q^{-2-\varepsilon}\; . $$ In terms of the convergents
$(p_k/q_k)_{k\ge 0}$ of $\theta$ with partial quotients
$(a_k)_{k\ge 1}$, this is equivalent to have, for all $\varepsilon
>0$ $$ a_{k+1}=\hbox{O}\, (q_k^{\varepsilon})\; , $$ which
explains our terminology.

We can reformulate (a) in terms of the lengths
$\lambda_\alpha^{(k)}$. It is convenient here to take as norm of a
matrix the sum of all coefficients (in absolute value; the
matrices that we consider here have nonnegative entries).

\vskip .3 truecm\noindent {\bf Proposition}{\it We have always,
for $k\ge 0$} $$ \hbox{Max}_{\alpha\in\cA}\lambda_\alpha^{(k)}\ge
\lambda^*\Vert
Q(k)\Vert^{-1}\ge\hbox{Min}_{\alpha\in\cA}\lambda_\alpha^{(k)}\; .
$$ {\it Condition (a) is equivalent to the following converse
estimate: for all $\varepsilon > 0$, there exists $C_\varepsilon
>0$ such that} $$\hbox{Max}_{\alpha\in\cA}\lambda_\alpha^{(k)}\le
C_\varepsilon \hbox{Min}_{\alpha\in\cA}\lambda_\alpha^{(k)}\Vert
Q(k)\Vert^{\varepsilon}\; . $$

\vskip .3 truecm \noindent \proof The first estimate follows from
$$
\lambda^*=\sum_{\alpha\in\cA}\lambda_\alpha^{(0)}=\sum_{\beta\in\cA}Q_\beta
(k)\lambda_\beta^{(k)}\; . $$ Assume (a) is satisfied. Let $l$ be
equal to $k+2d-3$ (if $d\ge 3$) or $k+2$ (if $d=2$) as in the
Lemma in 1.2.4. We have $$\Vert Q(k,l)\Vert\le C_\varepsilon '
\Vert Q(k)\Vert^{\varepsilon}\;  $$ for all $\varepsilon >0$ (with
an appropriate constant $C_\varepsilon '$).

This gives $$\hbox{Max}_{\alpha\in\cA}\lambda_\alpha^{(k)}\le
C_\varepsilon ' \Vert Q(k)\Vert^{\varepsilon}
\hbox{Max}_{\alpha\in\cA}\lambda_\alpha^{(l)}\; . $$

On the other hand, the Lemma 1.2.4 gives $$
\hbox{Min}_{\alpha\in\cA}\lambda_\alpha^{(k)}\ge
\hbox{Max}_{\alpha\in\cA}\lambda_\alpha^{(l)}\; , $$ giving the
required estimate. Assume now that the estimate of the Proposition
holds. We have always $$
\hbox{Max}_{\alpha\in\cA}\lambda_\alpha^{(k)}\ge d^{-1}\Vert
Z(k+1)\Vert \hbox{Min}_{\alpha\in\cA}\lambda_\alpha^{(k+1)}\; . $$
On the other hand, by definition of the $Z$ matrices, there exists
$\alpha_0\in\cA$ such that
$$\lambda_{\alpha_0}^{(k)}=\lambda_{\alpha_0}^{(k+1)}\; . $$

But we have $$\eqalign{\lambda_{\alpha_0}^{(k+1)} & \le
\hbox{Max}_{\alpha\in\cA}\lambda_\alpha^{(k+1)}\cr & \le
C_\varepsilon \hbox{Min}_{\alpha\in\cA} \lambda_\alpha^{(k+1)}
\Vert Q(k+1)\Vert^{\varepsilon} \cr &\le C_\varepsilon  d \Vert
Z(k+1)\Vert^{-1}  \hbox{Max}_{\alpha\in\cA}\lambda_\alpha^{(k)}
\Vert Q(k+1)\Vert^{\varepsilon} \cr &\le C_\varepsilon^2  d \Vert
Z(k+1)\Vert^{-1} \hbox{Min}_{\alpha\in\cA}\lambda_\alpha^{(k)}
\Vert Q(k)\Vert^{\varepsilon}\Vert Q(k+1)\Vert^{\varepsilon}\cr
&\le C_\varepsilon^2  d \Vert Z(k+1)\Vert^{-1}
\lambda_{\alpha_0}^{(k)} \Vert Q(k)\Vert^{\varepsilon}\Vert
Q(k+1)\Vert^{\varepsilon}\cr} $$ which implies $$ \Vert
Z(k+1)\Vert\le C_\varepsilon^2  d \Vert
Q(k)\Vert^{\varepsilon}\Vert Q(k+1)\Vert^{\varepsilon} $$ and
allows to conclude that (a) holds . \qed

\vskip .5 truecm\noindent {\bf Remark 1.}{\it Assume condition (a)
is satisfied. Set $k_0=2d-3$ if $d\ge 3$, $k_0=2$ if $d=2$.
Following the same lines that in the last Proposition, we see that
for any $\varepsilon >0$, there exists $C_\varepsilon >0$ such
that for $k\ge k_0$ we have $$\hbox{\rm Min}_{\alpha , \beta \in\cA}
Q_{\alpha\beta}(k)\ge C_\varepsilon^{-1}\Vert
Q(k)\Vert^{1-\varepsilon}\; . $$ On the other hand it is easy to
see that, even in the case of $3$ intervals this estimate {\it
does not} imply condition (a).}

\vskip .5 truecm\noindent {\bf Remark 2.}{ Boshernitzan has
defined ([Bo]) another condition which generalizes Roth condition
for irrational numbers. Namely, he asks that $T$ satisfies Keane's
condition and that the minimum distance $m_n$ between
discontinuity points of the $n$--th iterate $T^n$ of $T$  should
verify $$ m_n\ge {\gamma_\varepsilon^{-1}\over
n^{1+\varepsilon}}\; . $$ He proves that this condition has full
measure.

The relation between Boshernitzan's condition and condition (a)
above is however not clear.

\vskip .5 truecm \noindent {\bf 1.3.2 Spectral gap}

\vskip .3 truecm \noindent As soon as $k\ge 2d-3$ ($k\ge 2 $ if
$d=2$), all entries in the matrix $Q(k)$ are strictly positive. It
is therefore not unreasonable to expect that the positive cone is
more expanded by $Q(k)$ than the other directions, in the spirit
of Perron--Frobenius theorem.

However this is not automatic, as attested by the existence of
minimal non uniquely ergodic i.e.m.\ 's (an i.e.m.\ satisfying
Keane's condition is uniquely ergodic if and only if the image
under $Q(k)$ of the positive cone converges to a ray as
$k\rightarrow \infty$).

Our second condition ensures that this weird behaviour does not
occur.

For each $k\ge 0$, let $\Gamma^{(k)}$ be a copy of $\R^\cA$. One
should think of $\Gamma^{(k)}$ as the space of functions on
$\sqcup_{\alpha\in\cA} I_\alpha^{(k)}$ which are constant on each
$I_\alpha^{(k)}$. For $0\le k\le l$, let $S(k,l)$ be the linear
map from $\Gamma^{(k)}$ to $\Gamma^{(l)}$ whose matrix in the
canonical basis is $^tQ(k,l)$. This can be interpreted as a
special Birkhoff sum (see Section 2 below).

For $\varphi =(\varphi_\alpha )_{\alpha\in\cA}\in\Gamma^{(k)}$,
define $$ I_k(\varphi
)=\sum_{\alpha\in\cA}\lambda_\alpha^{(k)}\varphi_\alpha\; ; $$ we
have then $$ I_l(S(k,l)\varphi)=I_k(\varphi )\; . $$ Denote by
$\Gamma^{(k)}_*$ the kernel of the linear form $I_k$. We will ask
the following:
\item{(b)}{\it There exists $\theta >0, C>0$ such that, for all
$k\ge 0$, we have} $$ \Vert S(k)\mid_{\Gamma^{(0)}_*}\Vert\le
C\Vert S(k)\Vert^{1-\theta}=C\Vert Q(k)\Vert^{1-\theta}\; . $$

Observe that an i.e.m.\ satisfying Keane's condition and (b) must
be uniquely ergodic.

In appedix B we construct i.e.m.\ 's which satisfy condition (a)
but are {\it not} uniquely ergodic (see also [Ch]); therefore
condition (b) is not a consequence of condition (a).

However, if instead of condition (a) we consider the stronger
condition (reminding of bounded type irrational numbers):
\item{(\~ a)}{\it the sequence $Z(k)$ is bounded}

then condition (b) follows. Indeed, each $Q(k,k+2d-3)$ ($Q(k,k+2)$
when $d=2$) will contract by a definite factor $<1$ the Hilbert
metric of the projective positive cone.

\vskip .5 truecm \noindent {\bf 1.3.3 Coherence}

\vskip .3 truecm\noindent
To define our third condition, we
consider again the operators $S(k,l)\, :\,
\Gamma^{(k)}\rightarrow\Gamma^{(l)}$. Let $\Gamma^{(k)}_s$ be the
linear subspace of $\Gamma^{(k)}$ whose elements $v$ satisfy the
following: there exists $\sigma =\sigma (v)>0$, $C=C(v)>0$ such
that, for all $l\ge k$, one has $$\Vert S(k,l)v\Vert\le C\Vert
S(k,l)\Vert^{-\sigma}\Vert v\Vert\; . $$

We call $\Gamma^{(k)}_s$ the stable subspace of $\Gamma^{(k)}$.
Obviously, one has $\Gamma^{(k)}_s\subset\Gamma^{(k)}_*$. On the
other hand, $\Gamma^{(k)}_s$ is never reduced to $0$ because it
always contains the translation vector
$(\delta_\alpha^{(k)})_{\alpha\in\cA}$.

The operator $S(k,l)$ maps $\Gamma^{(k)}_s$ onto $\Gamma^{(l)}_s$.
Therefore we can define a quotient operator $$ S_\flat (k,l)\, :\,
\Gamma^{(k)}/\Gamma^{(k)}_s\rightarrow
\Gamma^{(l)}/\Gamma^{(l)}_s\; . $$ As we have quotiented out the
stable directions, it is not unreasonable to expect that the norm
of the inverse of $S_\flat (k,l)$ is not too large. This is what
our third condition is about:
\item{(c)}{\it for any $\varepsilon >0$, there exists
$C_\varepsilon >0$ such that, for all $l\ge k$, we have}
$$\eqalign{ \Vert [S_\flat (k,l)]^{-1}\Vert &\le
C_\varepsilon\Vert Q(l)\Vert^\varepsilon\; , \cr \Vert S
(k,l)\mid_{\Gamma^{(k)}_s}\Vert &\le C_\varepsilon\Vert
Q(l)\Vert^\varepsilon\; .\cr} $$

\vskip .3 truecm\noindent {\bf Remark.}{\it The second estimate in
(c) was wrongly omitted in [MMY].}

\vskip .5 truecm \noindent {\bf 1.3.4 Roth--type interval exchange
maps }

\vskip .3 truecm\noindent We say that an i.e.m.\ $T$ is of {\it
Roth type} if it satisfies Keane's condition and conditions (a),
(b), (c).

In the next Section, we will solve the cohomological equation for
i.e.m.\ 's of Roth type. In Section 4 we will prove the following

\vskip .3 truecm\noindent {\bf Theorem.}{\it Roth type interval
exchange maps form a subset of full measure.}

\vskip .3 truecm \noindent We also observe that if an i.e.m.\ $T$
satisfies Keane's condition, and its Rauzy--Veech continued
fraction is eventually periodic (meaning that the path $\gamma$ in
the Rauzy diagram is eventually periodic), then condition (\~ a),
(b) and (c) are automatically satisfied and therefore $T$ is of
Roth type.

\vskip 1. truecm \vfill\eject \noindent {\bf 2.The cohomological
equation}

\vskip .5 truecm \noindent {\bf 2.1 The Theorem of Gottschalk and
Hedlund }

\vskip .5 truecm \noindent {\bf 2.1.1 The statement}

\vskip .3 truecm \noindent We recall the following theorem of
Gottschalk and Hedlund. Let $X$ be a compact topological space,
$f$ a minimal homeomorphism of $X$ and $\psi$ a real valued
continuous function on $X$. Given $x_0\in X$ and $n\ge 1$ we
denote $S_n\psi(x_0)$ the Birkhoff sum $\sum_{j=0}^{n-1}\psi\circ
T^j(x_0)$. Suppose that there exists a point $x_0\in X$ and a
positive constant $C$ such that for all positive integer $n$ one
has $|S_n\psi(x_0)|\le C$. Then the cohomological equation $$
\varphi\circ f-\varphi=\psi $$ has a continuous solution
$\varphi$.

\vskip .5 truecm \noindent {\bf 2.1.2 Application to interval
exchange maps}

\vskip .3 truecm \noindent Let $T$ be an i.e.m.\ satisfying
Keane's condition. Then $T$ is minimal but not continuous.
However, the following well--known construction, reminiscent of
Denjoy counterexamples, allows to bypass this problem.

For $n\ge 0$, define $$\eqalign{ D_0(n) &= \{ T^{-n}(j_0(0,\alpha
))\, , \, \alpha\in\cA\, , \, \pi_0 (\alpha )>1\}\; , \cr D_1(n)
&= \{ T^{+n}(j_1(0,\alpha ))\, , \, \alpha\in\cA\, , \, \pi_1
(\alpha )>1\}\; . \cr} $$

It follows from the Keane property that these sets are disjoint
from each other and do not contain $0$.

Define an atomic measure $\mu$ by $$ \mu=\sum_{n\ge 0} \sum_{y\in
D_0(n)\sqcup D_1(n)} 2^{-n} \delta_y $$ and the increasing maps $
i^+,i^-\, :\, I\rightarrow\R$ by $$\eqalign{ i^-(y) &= y+\mu
([0,y))\; , \cr i^+(y) &= y+\mu ([0,y])\; . \cr}$$ We therefore
have $$\eqalign{i^+(y) &<i^-(y')\;\;\;\hbox{for}\, y<y'\; , \cr
i^+(y) &=i^-(y)\;\;\;\hbox{for}\, y\notin\sqcup_{n\ge
0}(D_0(n)\sqcup D_1(n))\; , \cr i^+(y) &=i^-(y)
+2^{-n}\;\;\;\hbox{for}\, y\in D_0(n)\sqcup D_1(n)\; . \cr}$$ We
also define $$
i^-(\lambda^*)=\lambda^*+4(d-1)=\lim_{y\nearrow\lambda^*}i^\pm(y)\;
, $$ and $$K=i^-(I)\cup
i^+(I)\cup\{i^-(\lambda^*)\}=\overline{i^-(I)}=\overline{i^+(I)}\;
.$$ As $T$ is minimal, $K$ is a Cantor set whose gaps are the
intervals $(i^-(y),i^+(y))$, $y\in\cup_{n\ge
0}\cup_{\varepsilon\in\{0,1\}}D_\varepsilon (n)$.

\vskip .3 truecm\noindent {\bf Proposition.}{\it There is a unique
continuous map $\hat T\, :\, K\rightarrow K$ such that $\hat
T\circ i^+=i^+\circ T$ on $I$. Moreover, $\hat T$ is a minimal
homeomorphism.}

\vskip .3 truecm\noindent The elementary proof is left to the
reader.

Let $\psi\, :\, I\rightarrow \R$ be a function which is continuous
on each $j_0(I_\alpha )$, with finite limits at the right
endpoints of each $j_0(I_\alpha )$. There is a unique continuous
function $\hat\psi\, ,\, K\rightarrow \R$ such that $\psi
(y)=\hat\psi\circ i^+ (y)$ for all $y\in I$. Assume that, for some
$x_0\in I$ the Birkhoff sums of $\psi$ for $T$ are bounded. Then
the same is true for the Birkhoff sums of $\hat\psi$ for $\hat T$
at the point $\hat{ x}_0=i^+(x_0)$. By the theorem of Gottschalk
and Hedlund, there is a continuous function $\hat\varphi\, ,\,
K\rightarrow \R$ satisfying $\hat\psi=\hat\varphi\circ\hat
T-\hat\varphi$. Define, for $y\in I$ $$\varphi
(y)=\hat\varphi\circ i^+ (y)\; . $$ In general, $\varphi$ is not
continuous. However it is bounded and satisfies $\varphi\circ
T-\varphi=\psi$. In the following, we will show that under
appropriate circumstances certain Birkhoff sums are bounded.

\vskip .5 truecm \noindent {\bf 2.2 Special Birkhoff sums }

\vskip .5 truecm \noindent {\bf 2.2.1 }

\vskip .3 truecm \noindent Let $T$ be an i.e.m.\ satisfying
Keane's condition. Denote by $T^{(k)}$ the i.e.m.\ obtained by the
accelerated Zorich algorithm (with $D=d-1$ in 1.2.4).

Let $\varphi\, :\, I^{(k)}\rightarrow\R$ be a function defined on
the domain $I^{(k)}$ of $T^{(k)}$. Let also $l\ge k$. For
$\beta\in\cA$, $x\in j_0(I_\beta^{(l)})$, the return time of $x$
into $I^{(l)}$ under iteration of $T^{(k)}$ is $Q_\beta (k,l)$.
Define a function $$ S(k,l)\varphi\, :\, I^{(l)}\rightarrow\R $$
by the formula $$S(k,l)(\varphi )(x)=\sum_{0\le i<Q_\beta
(k,l)}\varphi((T^{(k)})^i(x))\; , $$ for $x\in
j_0(I_\beta^{(l)})$. Observe that when $\varphi$ is constant on
each $ j_0(I_\beta^{(k)})$, the same is true of $S(k,l)(\varphi )$
in $I^{(l)}$ and the corresponding linear operator has $^tQ(k,l)$
as matrix in the canonical basis, as anticipated in 1.3.2.

We just write $S(k)$ for $S(0,k)$.

\vskip .5 truecm \noindent {\bf 2.2.2 Some elementary properties
of the operators $S(k,l)$}

\vskip .3 truecm \noindent
\item{2.2.2.1} For $m\ge l\ge k$ one has
$$ S(k,m)=S(l,m)\circ S(k,l)\; . $$
\item{2.2.2.2}
The operators $S(m,n)$ preserve all regularity classes which are
invariant by restriction, sum and translation.
\item{2.2.2.3}
If $\varphi$ is an integrable function on $I^{(k)}$,
 $$ \int_{I^{(k)}}\varphi(x)dx = \int_{I^{(l)}} (S(k,l)\varphi) (x)dx\; . $$
\item{2.2.2.4}
The operators $S(k,l)$ commute with taking derivatives.
\item{2.2.2.5} If the restriction of $\varphi$ to each
$j_0(I_\alpha^{(k)})$ is a polynomial of degree $\le\mu$, the
restriction of $S(k,l)\varphi$ to each $j_0(I_\beta^{(l)})$ is
also a polynomial of degree $\le\mu$. The case $\mu =0$ has
already been considered.
\item{2.2.2.6}Denote by $\hbox{BV}\, (\sqcup I_\alpha^{(k)})$ the space of
functions $\varphi$ on $I^{(k)}$ whose restriction to each
$j_0(I_\alpha^{(k)})$ has bounded variation and define $$
\hbox{Var}\, \varphi = \sum_\alpha \hbox{Var}\, \varphi
|_{j_0(I_\alpha^{(k)})}\; . $$ (We do not take into account the
discontinuities of $\varphi$ at the discontinuity points of
$T^{(k)}$). Then $S(k,l)$ sends $\hbox{BV}\, (\sqcup
I_\alpha^{(k)})$ into $\hbox{BV}\, (\sqcup I_\alpha^{(l)})$ and we
have $$ \hbox{Var}\, S(k,l)\varphi \le \hbox{Var}\, \varphi \; .
$$

\vskip .5 truecm\noindent {\bf 2.2.3 Reduction of Birkhoff sums to
special Birkhoff sums}

\vskip .3 truecm\noindent For diffeomorphisms of the circle with
irrational rotation number, when trying to estimate the Birkhoff
sums of some function, it is a standard trick to consider first
the ones associated to the denominators of the convergents of the
rotation number. We will do the same here.

Let $\varphi\, :\, I\rightarrow\R$ be a function, $x\in I$, and
$N>0$. We want to compute the Birkhoff sums $$S_N\varphi
(x)=\sum_{i=0}^{N-1} \varphi(T^i(x))$$ (with $T=T^{(0)}$).

We first replace $x$ by the point in the orbit $\{x,T(x), \ldots
,T^{N-1}(x)\}$ which is closest to $0$ and cut the Birkhoff sum
into two parts (one for $T$ and the other for $T^{-1}$). Let us
assume to keep notations simple that $x$ is actually closest to
the origin.

Let $k\ge 0$ be the largest integer such that at least one of the
points $T(x),\ldots ,T^{N-1}(x)$ belongs to $I^{(k)}$; because
$T^{(k)}$ is the first return map into $I^{(k)}$, these points are
precisely $T^{(k)}(x), (T^{(k)})^2(x),\ldots ,(T^{(k)})^{b(k)}(x)$
for some integer $b(k)>0$. Moreover, as none of these points
belongs to $I^{(k+1)}$ we must have
$$b(k)<\hbox{Max}_{\beta\in\cA}\sum_{\alpha\in\cA}Z_{\alpha\beta}(k+1)\;
, $$ the right hand term being the largest return time of
$T^{(k)}$ into $I^{(k+1)}$.

We set $x_k=x$, $x_{k-1}=(T^{(k)})^{b(k)}(x)$ and define
inductively $b(l)$ and $x_l$ for $0\le l\le k$.

The point $x_l$ has the property that it belongs to $I^{(l)}$ and
none of the points $T(x_l),\ldots ,T^N(x)$ belongs to $I^{(l+1)}$.
Those who belong to $I^{(l)}$ are $T^{(l)}(x_l),\ldots ,
(T^{(l)})^{b(l)}(x_l)$ for some integer $b(l)\ge 0$. We have $$
b(l)<\hbox{Max}_{\beta\in\cA}\sum_{\alpha\in\cA}Z_{\alpha\beta}(l+1)\;
, $$ We define $x_{l-1}=(T^{(l)})^{b(l)}(x_l)$. The process stops
when $x_l=T^N(x)$ (or $l=0$). From this construction it is obvious
that we have $$ S_N\varphi (x) = \sum_{l=0}^k\sum_{0\le i <b(l)}
S(l)\varphi((T^{(l)})^i(x_l))$$ which in particular implies, if
$\varphi$ is bounded: $$|S_N\varphi (x)|\le \sum_{l=0}^k\Vert
Z(l+1)\Vert\Vert S(l)\varphi\Vert_{L^\infty}\; , $$ where $\Vert
Z(l+1)\Vert = \hbox{Max}_{\beta\in\cA}\sum_{\alpha\in\cA}
Z_{\alpha\beta}(l+1)\; . $

In particular, if we are able to show that for some $\omega >0$ we
have $$\Vert S(l)\varphi\Vert_{L^\infty}\le C\Vert
Q(l)\Vert^{-\omega}\Vert\varphi\Vert\; , $$ and condition (a) in
1.3.1 is satisfied, then the Birkhoff sums of $\varphi$ will be
bounded.

\vskip .5 truecm \noindent {\bf 2.3 Estimates for functions of
bounded variation}

\vskip .5 truecm \noindent {\bf 2.3.1 }

\vskip .3 truecm \noindent Denote by $\hbox{BV}_*\, (\sqcup
I_\alpha^{(k)})$ the subspace of $\hbox{BV}\, (\sqcup
I_\alpha^{(k)})$ formed by the functions of mean $0$. The operator
$S(k,l)$ sends this subspace into $\hbox{BV}_*\, (\sqcup
I_\alpha^{(l)})$.

Let $\varphi\in\hbox{BV}_*\, (\sqcup I_\alpha^{(k)})$. We write $$
S(k,k+1)(\varphi )=\varphi_{k+1}+\chi_{k+1}\; , $$ with
$\chi_{k+1}\in\Gamma^{(k+1)}_*$ and $\varphi_{k+1}$ of mean zero
on {\it each} $j_0(I_\alpha^{(k+1)})$. Then we go on with:
$$S(j,j+1)(\varphi_j)=\varphi_{j+1}+\chi_{j+1}$$ with
$\chi_{j+1}\in\Gamma^{(j+1)}_*$ and $\varphi_{j+1}$ of mean zero
on {\it each} $j_0(I_\alpha^{(j+1)})$. We obtain, for $l>k$ $$
S(k,l)\varphi = S(l-1,l)\varphi_{l-1}+\sum_{k<j<l}S(j,l)\chi_j\; .
$$

\vskip .5 truecm \noindent {\bf 2.3.2 }

\vskip .3 truecm \noindent As $\varphi_j$ differs from
$S(k,j)\varphi$ by a function in $\Gamma^{(j)}_*$ and has mean
zero on each $j_0(I_\alpha^{(j)})$ we have (see 2.2.2.6) $$
\Vert\varphi_j\Vert_{L^\infty}\le \hbox{Var}\, \varphi_j\le
\hbox{Var}\, \varphi\; . $$ On the other hand, we have $$\eqalign{
\Vert\chi_j\Vert_{L^\infty} &\le \Vert
S(j-1,j)\varphi_{j-1}\Vert_{L^\infty}\cr &\le \Vert
Z(j)\Vert\Vert\varphi_{j-1}\Vert_{L^\infty}\cr}$$ (with
$\varphi_{j-1}=\varphi$ when $j=k+1$). We obtain therefore
$$\eqalign{ \Vert S(k,l)\varphi\Vert_{L^\infty} &\le \Vert
Z(k+1)\Vert\Vert
S(k+1,l)|_{\Gamma_*^{(k+1)}}\Vert\Vert\varphi_0\Vert\cr
&+\sum_{k<j<l}\Vert Z(j+1)\Vert\Vert
S(j+1,l)|_{\Gamma_*^{(j+1)}}\Vert\hbox{Var}\, \varphi\; . \cr}$$

\vskip .5 truecm \noindent {\bf 2.3.3 }

\vskip .3 truecm \noindent We now take $k=0$ and estimate the sum
$$\sum_{0\le j<l}\Vert Z(j+1)\Vert\Vert
S(j+1,l)|_{\Gamma_*^{(j+1)}}\Vert\; , $$ assuming that conditions
(a), (b) of Section 1.3 are satisfied.

On one side we have, by condition (a), for all $\varepsilon >0$:
$$ \Vert Z(j+1)\Vert\le C_\varepsilon\Vert Q(j)\Vert^\varepsilon\;
. $$ To estimate $\Vert S(j+1,l)|_{\Gamma_*^{(j+1)}}\Vert$ we
distinguish two cases. We assume condition (b) of 1.3.2, which
involves an exponent $1-\theta$ with $\theta >0$.
\item{i)} Assume first that $\Vert Q(j+1)\Vert\le\Vert
Q(l)\Vert^{\theta /d}$. As $Q(j+1)$ belongs to $\hbox{SL}\, (d,\Z
)$, we have $$\Vert (Q(j+1))^{-1}\Vert\le c\Vert Q(l)\Vert^{\theta
{d-1\over d}}\; . $$ Next we write $$ S(j+1,l)=S(0,l)\circ
(S(0,j+1))^{-1}\; , $$ and it follows from condition (b) that we
have $$\Vert S(j+1,l)|_{\Gamma_*^{(j+1)}}\Vert\le C\Vert
Q(l)\Vert^{1-{\theta\over d}}\; . $$
\item{ii)} Assume now that $\Vert Q(j+1)\Vert >\Vert
Q(l)\Vert^{\theta /d}$. If $l\le j+2d-2$ ($l\le j+3$ when $d=2$),
we just use $$\Vert S(j+1,l)|_{\Gamma_*^{(j+1)}}\Vert\le
C_\varepsilon\Vert Q(l)\Vert^\varepsilon $$ by condition (a). If
$l>j+2d-2$ ($l>j+3$), we write $$S(0,l)=S(j',l)S(j+1,j')S(0,j+1)$$
with $j'=j+2d-2$ ($j'=j+3$ when $d=2$). As the entries of
$Q(j+1,j')$ are positive integers we have $$\Vert Q(0,l)\Vert\ge
C\Vert Q(j',l)\Vert \Vert Q(0,j+1)\Vert\; , $$ which implies
$$\Vert Q(j',l)\Vert\le C\Vert Q(l)\Vert^{1-{\theta\over d}}\; .
$$ As we have also $$\Vert Q(j+1,j')\Vert\le C_\varepsilon\Vert
Q(j')\Vert^\varepsilon\; , $$ we obtain in this case that $$ \Vert
S(j+1,l)|_{\Gamma^{(j+1)}}\Vert\le C_\varepsilon\Vert
Q(l)\Vert^{1-{\theta \over d}+\varepsilon}\; . $$

Putting the two cases together and inserting this in the sum, we
obtain

\vskip .3 truecm\noindent {\bf Proposition}{\it For
$\varphi\in\hbox{BV}_* (\sqcup_{\alpha\in\cA}I_\alpha^{(0)})$,
$l\ge 0$, one has } $$\Vert S(l)\varphi\Vert_{L^\infty}\le C\Vert
Q(l)\Vert^{1-{\theta \over 2d}}\Vert\varphi\Vert_{BV}\; . $$

\vskip .3 truecm\noindent {\bf Remark} {\it In case i), the
estimate we got for $Q(j+1)^{-1}$ is far from optimal (it should
be of the order of $Q(j+1)$) but sufficient for our purposes.}

\vskip .5 truecm \noindent {\bf 2.4 Primitives of functions of
bounded variation}

\vskip .5 truecm \noindent {\bf 2.4.1 }

\vskip .3 truecm \noindent For $k\ge 0$, we will denote by
$\hbox{BV}^1\, (\sqcup_{\alpha\in\cA}I_\alpha^{(k)})$ the space of
functions $\varphi\, :\, I^{(k)}\rightarrow\R$ which are
absolutely continuous on each $j_0(I_\alpha^{(k)}))$ and whose
derivative on each $j_0(I_\alpha^{(k)}))$ is of bounded variation.
The condition that the mean value of the derivative is zero
defines an hyperspace $\hbox{BV}^1_*\,
(\sqcup_{\alpha\in\cA}I_\alpha^{(k)})$. We recall from 1.3.3 the
subspace $\Gamma_s^{(k)}$ of $\Gamma^{(k)}$. We will denote by
$\overline{\hbox{BV}}_*^1\, (\sqcup_{\alpha\in\cA}I_\alpha^{(k)})$
the quotient of $\hbox{BV}_*^1\,
(\sqcup_{\alpha\in\cA}I_\alpha^{(k)})$ by this finite dimensional
subspace.

Given $\varphi\in \hbox{BV}_*\,
(\sqcup_{\alpha\in\cA}I_\alpha^{(0)})$, we will find a primitive
$\Phi$ of $\varphi$ (given a priori by $d$ constants of
integration, one for each $I_\alpha^{(0)})$) for which the special
Birkhoff sums are small. The primitive $\Phi$ will actually be
uniquely determined $\hbox{mod}\, \Gamma_s^{(0)}$, i.e.\ in
$\overline{\hbox{BV}}_*^1\,
(\sqcup_{\alpha\in\cA}I_\alpha^{(0)})$.

\vskip .5 truecm \noindent {\bf 2.4.2 }

\vskip .3 truecm \noindent For any $\varphi\in \hbox{BV}_*\,
(\sqcup_{\alpha\in\cA}I_\alpha^{(k)})$, denote by
$P_0^{(k)}\varphi$ the class in $\overline{\hbox{BV}}_*^1\,
(\sqcup_{\alpha\in\cA}I_\alpha^{(k)})$ of the primitive of
$\varphi$ which has mean zero on each $j_0(I_\alpha^{(k)})$.

This is the most natural choice of primitive, but unfortunately
the special Birkhoff sums $S(k,l)$ do not commute with these
primitive operators, i.e.\ they do not preserve the condition to
be of mean value $0$ on each $j_0(I_\alpha^{(k)})$.

Therefore, we will modify $P_0^{(k)}$, considering $$
P^{(k)}=P_0^{(k)} +\Delta P^{(k)}\; , $$ where $$\Delta P^{(k)}\,
: \, \hbox{BV}_*\,
(\sqcup_{\alpha\in\cA}I_\alpha^{(k)})\rightarrow\Gamma^{(k)}/\Gamma^{(k)}_s\;
, $$ is a bounded linear operator. We want this new choice to be
equivariant: $$ S(k,l)\circ P^{(k)}=P^{(l)}\circ S(k,l)\; . $$
This leads to the following equation for $\Delta P^{(k)}$. Define
$$ \Lambda (k,l)= P_0^{(l)}\circ S(k,l)-S(k,l)\circ P_0^{(k)}\; .
$$ This is a bounded linear map from $\hbox{BV}_*\,
(\sqcup_{\alpha\in\cA}I_\alpha^{(k)})$ to
$\Gamma^{(l)}/\Gamma^{(l)}_s$. Then we should have $$ S_\flat
(k,l)\circ \Delta P^{(k)}- \Delta P^{(l)}\circ S(k,l)=\Lambda
(k,l)\; , \eqno(*) $$ where $S_\flat$ was defined in 1.3.3.

Equation $(*)$ has the formal solution $$ \Delta
P^{(k)}=\sum_{l>k}(S_\flat (k,l))^{-1}\circ \Lambda (l-1,l)\circ
S(k,l-1)\; , \eqno(**)$$ and we will check next that this defines
indeed the required primitive.

\vskip .5 truecm \noindent {\bf 2.4.3 Estimate for $\Lambda
(l-1,l)$. }

\vskip .3 truecm \noindent Let $\varphi\in \hbox{BV}_*\,
(\sqcup_{\alpha\in\cA}I_\alpha^{(l-1)})$. As $P_0^{(l-1)}\varphi$
has mean zero on each $j_0(I_\alpha^{(l-1)})$, we have $$ \Vert
P_0^{(l-1)}\varphi\Vert_{L^\infty}\le
(\hbox{Max}_{\alpha\in\cA}\lambda_\alpha^{(l-1)})\Vert\varphi\Vert_{L^\infty}\;
. $$ On the other hand we have $$ \eqalign{ \Vert
S(l-1,l)\varphi\Vert_{L^\infty} &\le \Vert
Z(l)\Vert\Vert\varphi\Vert_{L^\infty}\; , \cr \Vert
S(l-1,l)P_0^{(l-1)}\varphi\Vert_{L^\infty} &\le \Vert
Z(l)\Vert\Vert P_0^{(l-1)}\varphi\Vert_{L^\infty}\; . \cr}$$
Finally, we get $$\Vert
P_0^{(l)}S(l-1,l)\varphi\Vert_{L^\infty}\le
(\hbox{Max}_{\alpha\in\cA}\lambda_\alpha^{(l)})\Vert
Z(l)\Vert\Vert\varphi\Vert_{L^\infty}\;
, $$ which allows to conclude that $$\Vert\Lambda
(l-1,l)\varphi\Vert_{L^\infty}\le 2\Vert
Z(l)\Vert(\hbox{Max}_{\alpha\in\cA}\lambda_\alpha^{(l-1)})\Vert\varphi\Vert_{L^\infty}\;
. $$

\vskip .5 truecm \noindent {\bf 2.4.4  }

\vskip .3 truecm \noindent Assume now the three conditions (a),
(b) and (c) of 1.3. From 1.3.1, we get $$\eqalign{
\hbox{Max}_{\alpha\in\cA}\lambda_\alpha^{(l-1)} &\le C_\varepsilon
\Vert Q(l-1)\Vert^{\varepsilon -1}\; , \cr \Vert Z(l)\Vert &\le
C_\varepsilon \Vert Q(l-1)\Vert^{\varepsilon}\; , \cr} $$ and from
condition (c) that $$\Vert (S_\flat (0,l))^{-1}\Vert \le
C_\varepsilon \Vert Q(l)\Vert^{\varepsilon}\; . $$ On the other
hand, from the Proposition in 2.3, we obtain $$\Vert
S(0,l-1)\varphi\Vert_{L^\infty}\le C\Vert
Q(l-1)\Vert^{1-\theta/2d}\Vert\varphi\Vert_{BV}\; . $$ Therefore,
for $k=0$, the series (**) in 2.4.2 is converging and we obtain
$$\Vert \Delta P^{(0)}\varphi\Vert\le \left(\sum_{l>0}
C_\varepsilon '\Vert
Q(l)\Vert^{3\varepsilon-\theta/2d}\right)\Vert\varphi\Vert_{BV}\;
. $$ Indeed, we take $\varepsilon <\theta/6d$ and observe that it
follows from the Lemma in 1.2.4 that $\Vert Q(l)\Vert$ grows at
least exponentially fast. In the same way, as $T^{(k)}$ satisfies
also conditions (a), (b), (c) (with worse constants but the same
exponent $\theta$), the series (**) will converge for all $k\ge
0$. In this case, we prefer to estimate directly $\Delta
P^{(k)}S(0,k)\varphi$ for $\varphi\in\hbox{BV}_*\,
(\sqcup_{\alpha\in\cA}I_\alpha^{(0)})$. We have $$ \Vert\Delta
P^{(k)}S(0,k)\varphi\Vert\le \sum_{l>k}\Vert(S_\flat
(k,l))^{-1}\Vert\Vert\Lambda (l-1,l)\Vert\Vert
S(0,l-1)\varphi\Vert_{L^\infty}\; . $$ The above estimates now
give $$\eqalign{ \Vert\Delta P^{(k)}S(0,k)\varphi\Vert
&\le\left(\sum_{l>k}C_\varepsilon '\Vert Q(l)\Vert^{3\varepsilon
-\theta/2d}\right)\Vert\varphi\Vert_{BV}\cr &\le C \Vert
Q(k)\Vert^{-\theta /3d}\Vert\varphi\Vert_{BV}\; . \cr}$$

\vskip .5 truecm \noindent {\bf 2.4.5  Special Birkhoff sums for
$P^{(0)}\varphi$}

\vskip .3 truecm \noindent Let $\varphi\in \hbox{BV}_*\,
(\sqcup_{\alpha\in\cA}I_\alpha^{(0)})$, $\Phi\in \hbox{BV}_*^1\,
(\sqcup_{\alpha\in\cA}I_\alpha^{(0)})$ such that the class mod
$\Gamma_s^{(0)}$ of $\Phi$ is $P^{(0)}\varphi$. The class mod
$\Gamma_s^{(k)}$ of $S(k)\Phi$ is $P^{(k)}S(k)\varphi$ by
construction.

From the definition of $P_0^{(k)}$ and 2.3, we have $$\Vert
P_0^{(k)}S(k)\varphi\Vert_{L^\infty}\le\left(
\hbox{Max}_{\alpha\in\cA}\lambda_\alpha^{(k)}\right)\Vert
Q(k)\Vert^{1-\theta /d}\Vert\varphi\Vert_{BV}\; , $$ with
$\hbox{Max}_{\alpha\in\cA}\lambda_\alpha^{(k)}\le
C_\varepsilon\Vert Q(k)\Vert^{\varepsilon -1}$ by condition (a).
Joining this with the estimate for $\Delta P^{(k)}$ above, we
obtain $$\eqalign{\Vert S(k)P^{(0)}\varphi\Vert &= \Vert
P^{(k)}S(k)\varphi\Vert \cr &\le C\Vert
Q(k)\Vert^{-\theta/3d}\Vert\varphi\Vert_{BV}\; . \cr}$$ By
definition of a quotient norm, this means that we may write in
$\hbox{BV}_*^1\, ($ $\sqcup_{\alpha\in\cA}$ $I_\alpha^{(k)})$ : $$
S(k)\Phi = \Phi_k+\chi_k\; , $$ with $\chi_k\in\Gamma_s^{(k)}$ and
$$\Vert\Phi_k\Vert_{L^\infty}\le C\Vert Q(k)\Vert^{-\theta
/3d}\Vert\varphi\Vert_{BV}\; . $$ we have then $$\eqalign{
\chi_{k+1} &= S(k,k+1)\chi_k+ S(k,k+1)\Phi_k - \Phi_{k+1} \cr & :=
S(k,k+1) \chi_k +\Delta \chi_{k+1}\; , \cr}$$ with
$\Vert\Delta\chi_{k+1}\Vert\le C\Vert Q(k+1)\Vert^{-\theta
/4d}\Vert\varphi\Vert_{BV}$ (using once more condition (a)). Then
$$S(k)\Phi =\Phi_k +\sum_{j\le k} S(j,k)\Delta \chi_j\; . $$ In
the sum, we separate two cases. Recall that there exists $\sigma
>0$, $C>0$, such that $$\left\Vert
S(j)|_{\Gamma_s^{(0)}}\right\Vert\le C\Vert Q(j)\Vert^{-\sigma}\;
, $$ for all $j\ge 0$. If $\Vert Q(j)\Vert\le\Vert
Q(k)\Vert^{\sigma /d}$, we write $S(j,k)=S(k)\circ (S(j))^{-1}$
and get $$\eqalign{ \Vert S(j,k)\Delta\chi_j\Vert &\le C\Vert
Q(k)\Vert^{-\sigma}\Vert (S(j))^{-1}\Delta\chi_j\Vert \cr &\le
C\Vert Q(k)\Vert^{-\sigma} \Vert
Q(j)\Vert^{d-1}\Vert\Delta\chi_j\Vert \cr &\le C\Vert
Q(k)\Vert^{-\sigma /d}\Vert\varphi\Vert_{BV}\; . \cr}$$ In case
$\Vert Q(j)\Vert >\Vert Q(k)\Vert^{\sigma /d}$, we use the second
estimate in condition (c) to get $$\eqalign{ \Vert
S(j,k)\Delta\chi_j\Vert &\le C_\varepsilon \Vert
Q(k)\Vert^{\varepsilon }\Vert\Delta\chi_j\Vert \cr &\le
C_\varepsilon \Vert Q(k)\Vert^{\varepsilon -\theta\sigma /4d^2}
\Vert\varphi\Vert_{BV}\; . \cr}$$

We have thus proved the

\vskip .3 truecm\noindent {\bf Theorem.} {\it Let $T$ be an
i.e.m.\ of Roth type. There exists $\omega >0$, depending only on
$\sigma$ and $\theta$ in (b), (c), such that the special Birkhoff
sums $S(k)\Phi$ satisfy: } $$ \Vert S(k)\Phi\Vert_{L^\infty}\le
C\Vert Q(k)\Vert^{-\omega}\Vert\varphi\Vert_{BV}\; . $$

\vskip .3 truecm\noindent {\bf Corollary.} {\it Let $T$ be an
i.e.m.\ of Roth type, $\varphi\in \hbox{BV}_*\,
(\sqcup_{\alpha\in\cA}I_\alpha^{(0)})$. For any primitive $\Phi$
of $\varphi$ whose class lie in $P^{(0)}\varphi$, we can solve the
cohomological equation } $$ \Psi\circ T-\Psi = \Phi $$ {\it with a
bounded solution $\Psi$.}

\vskip .3 truecm\noindent \proof This follows from the Theorem,
taking into account the remarks at the end of Section 2.1 and
2.2.3.

\vfill\eject \noindent {\bf 3. Suspensions of interval exchange
maps}

\vskip .5 truecm \noindent We first recall, basically to fix
notations, how to suspend i.e.m.\ 's in order to get a Riemann
surface with an holomorphic $1$--form. The basic reference is
[V1].

\vskip .5 truecm \noindent {\bf 3.1 Suspension data }

\vskip .5 truecm \noindent Let $(\cA ,\pi_0 ,\pi_1)$ be admissible
combinatorial data, and let $T$ be an i.e.m.\ of this
combinatorial type, determined by length data
$(\lambda_\alpha)_{\alpha\in\cA}$. We will construct a Riemann
surface with a flow which can be considered as a suspension of
$T$. In order to do this, we need data which we call {\it
suspension data}. We will identify $\R^2$ with $\C$. Consider a
family $\tau =(\tau_\alpha)_{\alpha\in\cA}\in\R^\cA$. To this
family we associate $$\eqalign{ \zeta_\alpha &= \lambda_\alpha
+i\tau_\alpha \; , \; \; \alpha\in\cA\cr \xi_\alpha^\varepsilon &=
\sum_{\pi_\varepsilon \beta\le\pi_\varepsilon\alpha}\zeta_\beta\;
, \; \; \alpha\in\cA\; , \; \; \varepsilon\in\{0,1\}\; . \cr}$$ We
always have $\xi_{\alpha_0}^0=\xi_{\alpha_1}^1$, where as before
$\pi_\varepsilon (\alpha_\varepsilon )=d$. We say that $\tau$
defines suspension data if the following inequalities hold: $$
\eqalign{ \IM \xi_\alpha^0  & > 0 \; \; \hbox{for all}\,
\alpha\in\cA\, , \; \alpha\not=\alpha_0\; , \cr \IM \xi_\alpha^1 &
< 0 \; \; \hbox{for all}\, \alpha\in\cA\, , \;
\alpha\not=\alpha_1\; . \cr}$$ We also set $$ \theta_\alpha =
\xi_\alpha^1-\xi_\alpha^0\; , \alpha\in\cA\; . $$ We then have
$$\eqalign{ \theta &=\Omega\zeta\; , \cr \RE\theta &=\delta\; ,
\cr}$$ and define $$h=-\IM\theta =-\Omega\tau\; . $$ One has
$h_\alpha >0$ for all $\alpha\in\cA$, because of the formula
$$\theta_\alpha =
(\xi_\alpha^1-\zeta_\alpha)-(\xi_\alpha^0-\zeta_\alpha)\; . $$ One
has also $$\IM\xi_{\alpha_0}^0 = \IM\xi_{\alpha_1}^1\in
[-h_{\alpha_1},h_{\alpha_0}]\;  . $$

\vskip .5 truecm \noindent {\bf 3.2 Construction of a Riemann
surface }

\vskip .5 truecm \noindent Let $(\cA ,\pi_0 ,\pi_1)$ and
$(\zeta_\alpha=\lambda_\alpha +i\tau_\alpha )_{\alpha\in\cA}$ as
above. For $\alpha\in\cA$, consider the rectangles in $\C=\R^2$:
$$\eqalign{ R_\alpha^0 &= (\RE\xi_\alpha^0-\lambda_\alpha ,
\RE\xi_\alpha^0)\times [0,h_\alpha ]\; , \cr R_\alpha^1 &=
(\RE\xi_\alpha^1-\lambda_\alpha , \RE\xi_\alpha^1)\times
[-h_\alpha , 0 ]\; , \cr}$$ and the segments $$ \eqalign{
S_\alpha^0 &= \{\RE\xi_\alpha^0\}\times [0,\IM\xi_\alpha^0)\; , \;
\alpha\not=\alpha_0\; , \cr S_\alpha^1 &=
\{\RE\xi_\alpha^1\}\times (\IM\xi_\alpha^1, 0]\; , \;
\alpha\not=\alpha_1\; . \cr}$$ Let also
$S^0_{\alpha_0}=S^1_{\alpha_1}$ be the half--open vertical segment
$[\lambda^*,\xi_{\alpha_0}^0)=[\lambda^*,\xi_{\alpha_1}^1)$.

Define then $$R_\zeta =
\cup_{\varepsilon\in\{0,1\}}\cup_{\alpha\in\cA}R_\alpha^\varepsilon
\cup_{\varepsilon\in\{0,1\}}\cup_{\alpha\in\cA}S_\alpha^\varepsilon\;
. $$ The translation by $\theta_\alpha$ sends $R_\alpha^0$ onto
$R_\alpha^1$. If $\xi_{\alpha_0}^0=\xi_{\alpha_1}^1=0$,
$S^0_{\alpha_0}=S^1_{\alpha_1}$ is empty, $\xi_{\alpha_1}^0$ is
the top right corner of $R_{\alpha_1}^0$ and $\xi_{\alpha_0}^1$ is
the bottom right corner of $R_{\alpha_0}^1$. If
$\xi_{\alpha_0}^0=\xi_{\alpha_1}^1>0$, the translation by
$\theta_{\alpha_1}$ sends the top part
$\tilde{S}^0_{\alpha_1}=\{\RE\xi_{\alpha_1}^0\}\times
[h_{\alpha_1},\IM\xi_{\alpha_1}^0 )$ of $S^0_{\alpha_1}$ onto
$S^1_{\alpha_1}$. If $\xi_{\alpha_0}^0=\xi_{\alpha_1}^1<0$, the
translation by $\theta_{\alpha_0}$ sends $S^0_{\alpha_0}$ onto the
bottom part $\tilde{S}^1_{\alpha_0}=\{\RE\xi_{\alpha_0}^1\}\times
( \IM\xi_{\alpha_0}^1 , -h_{\alpha_0}]$ of $S^1_{\alpha_0}$.

We use these translations to identify in $R_\zeta$ each
$R^0_\alpha$ to each $R^1_\alpha$, and
$S^0_{\alpha_0}=S^1_{\alpha_1}$ (if non empty) to either
$\tilde{S}^0_{\alpha_1}$ or $\tilde{S}^1_{\alpha_0}$.

Denote by $M_\zeta^*$ the topological space obtained from
$R_\zeta$ by these identifications.

Observe that $M_\zeta^*$ inherits from $\C$ the structure of a
Riemann surface, and also a nowhere vanishing holomorphic
$1$--form $\omega$ (given by dz) and a vertical vector field
(given by ${\partial\over\partial y}$).

\vskip .5 truecm \noindent {\bf 3.3 Compactification of
$M_\zeta^*$ }

\vskip .5 truecm \noindent Let $\overline{\cA}$ be the set with
$2d-2$ elements of pairs $(\alpha ,L)$ and $(\alpha ,R)$, except
that we identify $(\alpha_0,R)=(\alpha_1,R)$ and
$(\alpha_0',L)=(\alpha_1',L)$, where $\pi_\varepsilon
(\alpha_\varepsilon )=d$, $\pi_\varepsilon (\alpha_\varepsilon '
)=1$.

Let $\sigma$ be the permutation of $\overline{\cA}$ defined by
$$\eqalign{ \sigma (\alpha , R) &= (\beta_0 ,L)\; , \cr \sigma
(\alpha , L) &= (\beta_1 ,R)\; , \cr}$$ with $\pi_0
(\beta_0)=\pi_0(\alpha )+1$, $\pi_1 (\beta_1)= \pi_1 (\alpha )-1$;
in particular, we have $$\eqalign{\sigma (\alpha_0 , R) &=
(\pi_0^{-1}(\pi_0(\alpha_1 )+1) ,L)\; , \cr \sigma (\alpha_1' , L)
&= (\pi_1^{-1}(\pi_1 (\alpha_0' )-1) ,R)\; . \cr}$$ The
permutation describes which half planes are met when one winds
around an end of $M_\zeta^*$. Denote by $\Sigma$ the set of cycles
of $\sigma$. To each $C\in\Sigma$ is associated in a one--to--one
correspondance an end $q_C$ of $M_\zeta^*$. From the local
structure around $q_C$, it is clear that the compactification
$M_\zeta = M_\zeta^*\cup_{C\in \Sigma}\{q_C\}$ will be a compact
Riemann surface, with the set of marked points $\cup_{C\in
\Sigma}\{q_C\}=M_\zeta \setminus M_\zeta^*$ in canonical
correspondance with $\Sigma$. Moreover, the $1$--form $\omega$
extends to a holomorphic $1$--form on $M_\zeta$; the length of a
cycle $C$ is an even number $2n_C$; the corresponding marked point
$q_C$ is a zero of $\omega$ of order $n_C-1$.

Let $\nu=\hbox{card}\, \Sigma$, and let $g$ be the genus of
$M_\zeta$. We have $$\eqalign{ d-1 &=\sum_{C\in\Sigma} n_C\; , \cr
2g-2 &=\sum_{C\in\Sigma} (n_C-1)\; , \cr}$$ hence $$d=2g+\nu -1\;
. $$

\vskip .3 truecm\noindent {\bf Example} Suppose that $\pi_0$,
$\pi_1$ satisfy $$\pi_0(\alpha)+\pi_1(\alpha )=d+1\;\; , \;
\hbox{for all}\, \alpha\in\cA\; . $$ If $d$ is even, there is only
$1$ cycle; we have $d=2g$ and the only zero of $\omega$ has order
$2g-2$. If $d$ is odd, there are two cycles of equal length $d-1$;
we have $d=2g+1$, and each of the two zeros of $\omega$ has order
$g-1$.

\vskip .3 truecm\noindent The vertical vector field on $M_\zeta^*$
does not extend (continuously) to $M_\zeta$ when $g>1$, unless one
slows it near the marked points (which we will not do here).
Nevertheless, it can be considered as a suspension of $T$:
starting from a point $(x,0)$ on the bottom side of $R_\alpha^0$,
one flows up till reaching the top side where the point
$(x,h_\alpha)$ is identified with the point  $(x+\delta_\alpha
,0)=(T(x),0)$ in the top side of $R_\alpha^1$. The return time is
$h_\alpha$. The vector field is not complete, as some orbits reach
marked points in finite time.

\vskip .5 truecm \noindent {\bf 3.4 The cohomological equation for
higher smoothness }

\vskip .5 truecm \noindent {\bf 3.4.1}

\vskip .3 truecm\noindent In this section, we will relate the
(discrete) cohomological equation for i.e.m.\ 's to the
(continuous) cohomological equation for the vertical vector field
on $M_\zeta$; this equation is $$\tilde{\Phi} =
{\partial\over\partial y}\tilde{\Psi}\; , $$ where now
$\tilde{\Phi}$, $\tilde{\Psi}$ are functions on $M_\zeta$. This
allows to compare our results with the pioneering work of Forni
([Fo1]). We will always assume, as he does, that $\tilde{\Phi}$
vanishes in the neighborhood of the marked points of $M_\zeta$.

Considering the cohomological equation on the surface leads
naturally to some regularity assumptions on the interval. Because
the datum $\Phi$ and the solution $\Psi$ are not related to the
corresponding functions $\tilde{\Phi}$, $\tilde{\Psi}$ on the
surface in the same way ($\Psi$ is a restriction of $\tilde{\Psi}$
to a segment, while $\Phi$ is an integral), the regularity that we
introduce for $\Phi$ and $\Psi$ are not of the same kind (even
taking the loss of derivatives into account).

\vskip .5 truecm \noindent {\bf 3.4.2}

\vskip .3 truecm\noindent For each integer $r\ge 1$, we introduce
the space $\hbox{BV}^r_*(I)$ of functions $\Phi\, :\,
I\rightarrow\R$ such that
\item{$\bullet$} for each $\alpha\in\cA$, $\Phi$ is of class
$\cC^{r-1}$ on $j_0(I_\alpha )$, $D^{r-1}\Phi$ is absolutely
continuous on $j_0(I_\alpha )$ and $D^r\Phi$ is of bounded
variation on $j_0(I_\alpha )$;
\item{$\bullet$} each function $D^l\Phi$, for $0<l\le r$, has mean
value $0$ in $I$.

\vskip .3 truecm\noindent {\bf Remark.}{\it As before, we allow
discontinuities at the discontinuities of $T$. Observe however
that the mean value condition implies that the sum of the jumps of
$D^l\Phi$ ($0\le l<r$) over the discontinuities of $T$ (including
the endpoints of $I$) is zero.}

\vskip .3 truecm\noindent We will indicate below why the mean
value condition is natural.

On the other hand, we will look for solutions in the space
$\cC^{r-2+Lip}(I)$ of functions $\Psi$ which are $\cC^{r-2}$ {\it
on all of $I$}, the derivative of order $r-2$ $D^{r-2}\Psi$ being
Lipschitz on $I$. For $r=1$, this is just the space of bounded
functions on $I$. Observe that, as soon as $r\ge 2$, we do not
allow discontinuities.

\vskip .5 truecm \noindent {\bf 3.4.3}

\vskip .3 truecm\noindent For $T$ an i.e.m.\ of Roth type, denote
by $\Gamma_T=\Gamma_T^{(0)}$ the space of functions
$\chi\in\Gamma$ (constant on each $j_0(I_\alpha )$) which can be
written as $$\chi=\psi-\psi\circ T$$ with bounded $\psi$. This is
a linear subspace of $\Gamma$ which is contained in $\Gamma_*$ and
contains $\Gamma_s$. We can rephrase our main theorem by saying
that there is a well--defined {\it obstruction map} $$
\hbox{BV}^1_*(I)\rightarrow \Gamma /\Gamma_T\; $$ which associates
to $\Phi$ the function in $\Gamma$ we must subtract from $\Phi$ in
order to be able to solve the cohomological equation. We recognize
(some of ) Forni distribution conditions, by choosing a basis in
the finite--dimensional space $ \Gamma /\Gamma_T$. The number of
conditions is just the codimension of $\Gamma_T$, as the
restriction of the obstruction map to $\Gamma$ is just the quotient map
and thus
the obstruction map is onto.

\vskip .5 truecm \noindent {\bf 3.4.4}

\vskip .3 truecm\noindent Let now $r\ge 1$, $\Phi\in
\hbox{BV}^r_*(I)$, and let us try to solve (under finitely many
linear conditions on $\Phi$) the cohomological equation $$\Phi =
\Psi-\Psi\circ T\; , $$ with $\Psi\in\cC^{r-2+Lip}(I)$. We assume
that the i.e.m.\ is of Roth type. Consider the $rd$--dimensional
space $\Gamma (r)$ of functions $\chi$ on $I$ whose restrictions
to each $j_0(I_\alpha )$ are polynomials of degree $<r$. For
$r=1$, this is our previous space $\Gamma$. Consider also $$
\Gamma_* (r)=\Gamma (r)\cap\hbox{BV}^r_*(I)$$ which has
codimension $(r-1)$ in $\Gamma (r)$. We first describe the
subspace $\Gamma_T (r)$ of $\Gamma_* (r)$ of functions $\chi$
which can be written as $$\chi = \psi\circ T-\psi\; , $$ with
$\psi\in\cC^{r-2+Lip}(I)$.

\vskip .5 truecm \noindent {\bf Lemma} {\it For $r\ge 1$, the map
$\chi\mapsto D\chi$ from  $\Gamma (r+1)$ to $\Gamma (r)$ sends
$\Gamma_* (r+1)$ to $\Gamma_* (r)$ and $\Gamma_T (r+1)$ to
$\Gamma_T (r)$. The kernel, i.e.\ the intersection
$\Gamma\cap\Gamma_T (r+1)$, is equal to $\R\delta$; we have thus}
$$ \hbox{dim}\,\Gamma_T (r)=\hbox{dim}\,\Gamma_T+(r-1)\; .$$

\vskip .3 truecm \noindent\proof It is clear that $\chi\mapsto
D\chi$ sends $\Gamma_* (r+1)$ to $\Gamma_* (r)$ and $\Gamma_T
(r+1)$ to $\Gamma_T (r)$. If $\psi_0(x) \equiv x$, then
$\psi_0\circ T(x)-\psi_0(x)=\delta_\alpha$ for $x\in j_0(I_\alpha
)$ hence $\R\delta \subset \Gamma\cap\Gamma_T (r)$ for all $r\ge
1$. Conversely, if $\chi\in\Gamma\cap\Gamma_T (r)$, write
$\chi=\psi\circ T-\psi$ with $\psi\in\hbox{Lip}\, (I)$. Taking
derivatives, $D\psi$ is $T$--invariant, hence constant as $T$ is
ergodic. Therefore $\chi\in\R\delta$. \qed

\vskip .5 truecm \noindent {\bf Theorem} {\it Let $r\ge 1$. For
any $\Phi\in\hbox{BV}^r_*(I)$, one can find $\chi\in\Gamma_* (r)$,
$\psi\in\cC^{r-2+Lip}(I)$ such that $$\Phi=\chi+\Psi\circ T-\Psi\;
.$$ In other terms, the map $$\eqalign{ \Phi & \mapsto \chi\cr
\hbox{BV}^r_*(I) & \rightarrow\Gamma_* (r)/\Gamma_T (r)\; , \cr}$$
is the obstruction map associated with the cohomological equation
with the prescribed regularities.}

\vskip .3 truecm \noindent\proof By induction on $r$, the case
$r=1$ being our main theorem. Assume
$\Phi\in\hbox{BV}^{r+1}_*(I)$. Then $D\Phi\in\hbox{BV}^r_*(I)$. By
the induction hypothesis, one can write $$D\Phi=\chi_1+\Psi_1\circ
T-\Psi_1\; , $$ with $\chi_1\in \Gamma_* (r)$ and
$\Psi_1\in\cC^{r-2+Lip}(I)$. Let $\Psi$ be a primitive of
$\Psi_1$, $\chi_0$ be a primitive of $\chi_1$. Then
$\psi\in\cC^{r-1+Lip}(I)$. As $D\Phi$ has mean value $0$, $\chi_1$
has also mean value $0$ and $\chi_0\in \Gamma_* (r+1)$. The
difference $\chi_0'=\Phi-\chi_0- \Psi\circ T+\Psi$ belongs to
$\Gamma$ and we take $\chi=\chi_0+\chi_0'$. \qed

\vskip .5 truecm \noindent {\bf 3.4.5}

\vskip .3 truecm\noindent We explain now why the regularities for
$\Phi$, $\Psi$ are ``natural''.

Let $\zeta=(\zeta_\alpha )_{\alpha\in\cA}$ be suspension data, and
let $M_\zeta$ be the surface constructed from these data as in
3.2.

Let $\tilde{\Phi}$ be a continuous function on $M_\zeta$. With the
notations of 3.1, we define, for $\alpha\in\cA$: $$ \eqalign{
\cI_\alpha^0 &=
\int_0^{\IM\xi_\alpha^0}\tilde{\Phi}(\RE\xi_\alpha^0,y)dy\; , \cr
\cI_\alpha^1 &=
\int^0_{\IM\xi_\alpha^1}\tilde{\Phi}(\RE\xi_\alpha^1,y)dy\; ;
\cr}$$ for $\alpha\in\cA$, $x\in j_0(I_\alpha )$, we also set
$$\Phi (x)=\int_0^{h_\alpha}\tilde{\Phi}(x,y)dy\;.$$ Observe that
we have $$\eqalign{ \Phi ((\RE\xi_\alpha^0)^-) &=
\cI_\alpha^0+\cI_\alpha^1\; , \cr \Phi
(\RE\xi_\alpha^0-\lambda_\alpha) &=
\cI_{\beta_0}^0+\cI_{\beta_1}^1\; , \cr}$$ where
$\pi_0(\beta_0)+1=\pi_0(\alpha)$,
$\pi_1(\beta_1)+1=\pi_1(\alpha)$, except if $\pi_0(\alpha)=1$
(respectively $\pi_1(\alpha)=1$) when $\cI_{\beta_0}^0$ (resp.\
$\cI_{\beta_1}^1$) is declared to be $0$.

From these formulas and $\cI_{\alpha_0}^0+\cI_{\alpha_1}^1=0$
(with $\pi_\varepsilon (\alpha_\varepsilon) =d $ as usual ), we
obtain $$ \sum_{\alpha\in\cA}\Phi ((\RE\xi_\alpha^0)^-)=
\sum_{\alpha\in\cA}\Phi (\RE\xi_\alpha^0-\lambda_\alpha )\; , $$
which means that the derivative of $\Phi$ (when it exists) has
mean value $0$. This explains the conditions defining
$\hbox{BV}^{r}_*(I)$. On the other hand, if $\tilde{\Psi}$ is a
function on $M_\zeta$ satisfying $${\partial\over\partial
y}\tilde{\Psi}=\tilde{\Phi}$$ and we define $$\Psi (x) =
\tilde{\Psi}(x,0) $$ then we will have $$\Psi\circ T-\Psi=\Phi\; .
$$

\vfill\eject \noindent {\bf 4. Proof of full measure for Roth type}

\vskip .5 truecm \noindent We will first recall the construction
of the finite measure, absolutely continuous w.r.t.\ Lebesgue
measure, which is invariant under the Zorich algorithm
(normalized).

\vskip .5 truecm \noindent {\bf 4.1 The basic operation of the
algorithm for suspensions}

\vskip .5 truecm \noindent Let $(\cA ,\pi_0,\pi_1)$ and
$(\zeta_\alpha = \lambda_\alpha+i\tau_\alpha )_{\alpha\in\cA}$ be
as above. Construct $R_\zeta$, $M_\zeta$ as in 3.2 and 3.3. With
$\pi_\varepsilon (\alpha_\varepsilon )=d$ as above, assume that
$$\lambda_{\alpha_0}\not=\lambda_{\alpha_1}\; . $$ Then the
formula $\lambda_{\alpha_\varepsilon}=\hbox{Max}\,
(\lambda_{\alpha_0} , \lambda_{\alpha_1})$ defines uniquely $
\varepsilon\in\{0,1\}$ and determines uniquely the basic step of
the continued fraction algorithm; this step produces new
combinatorial data $(\cA ,\hat\pi_0,\hat\pi_1)$ and length data
$(\hat\lambda_\alpha )_{\alpha\in\cA}$ given by $$ \eqalign{
\hat\lambda_\alpha &= \lambda_\alpha\; , \;\;
\alpha\not=\alpha_\varepsilon\cr \hat\lambda_{\alpha_\varepsilon}
&=\lambda_{\alpha_\varepsilon}-\lambda_{\alpha_{1-\varepsilon}}\;
. \cr}$$

For suspension data, we just define in the same way $$ \eqalign{
\hat\zeta_\alpha &= \zeta_\alpha\; , \;\;
\alpha\not=\alpha_\varepsilon\cr \hat\zeta_{\alpha_\varepsilon}
&=\zeta_{\alpha_\varepsilon}-\zeta_{\alpha_{1-\varepsilon}}\; .
\cr}$$ This has a nice representation in terms of the
corresponding regions $R_\zeta$, $R_{\hat\zeta}$. One cuts from
$R_\zeta$ the part where $x\ge \hat
\lambda^*=\lambda^*-\lambda_{\alpha_\varepsilon}$: it is made of
$R_{\alpha_{1-\varepsilon}}^{1-\varepsilon}$ and a right part of
$R_{\alpha_{\varepsilon}}^{\varepsilon}$. We glue back
$R_{\alpha_{1-\varepsilon}}^{1-\varepsilon}$ to the free
horizontal side of $R_{\alpha_{\varepsilon}}^{1-\varepsilon}$, and
the right part of $R_{\alpha_{\varepsilon}}^{\varepsilon}$ to
$R_{\alpha_{1-\varepsilon}}^{\varepsilon}$: see figure 4.


It is easy to check that the new suspension data satisfy the
inequalities required in 3.1; if for instance $\varepsilon=0$, one
has $$\hat\xi_\alpha^0=\xi_\alpha^0\;\; , \; \alpha\not=\alpha_0$$
with $\hat\pi_0=\pi_0$ on one hand and $$\eqalign{
\hat\xi_\alpha^1 &=\xi_\alpha^1\;\; , \; \alpha\not=\alpha_0,
\alpha_1\cr \hat\xi_{\alpha_{1}}^{1} &= \xi_{\alpha_{0}}^{1}\;\;
,\cr \hat\xi_{\alpha_{0}}^{1} &=
\xi_{\alpha_{0}}^{1}-\zeta_{\alpha_{1}}\;\; . \cr}$$ The last
formula gives $$\eqalign{ - \hat\xi_{\alpha_{0}}^{1} &=
\zeta_{\alpha_{1}}-\xi_{\alpha_{0}}^{1} \cr &=
\zeta_{\alpha_{1}}-\xi_{\alpha_{0}}^{0}-\theta_{\alpha_{0}}\cr &=
\zeta_{\alpha_{1}}-\xi_{\alpha_{1}}^{1}-\theta_{\alpha_{0}}\cr &=
-\xi_{\tilde\alpha_{1}}^{1}-\theta_{\alpha_{0}}\; , \cr}$$ with
$\pi_1(\tilde\alpha_{1})=d-1$. We therefore have $$-\IM
\hat\xi_{\alpha_{0}}^{1} = - \IM
\xi_{\tilde\alpha_{1}}^{1}+h_{\alpha_{0}}>0\; . $$ We also see
that (still with $\varepsilon =0$), if $\hat\alpha_1\in\cA$ is
such that $\hat\pi_1(\hat\alpha_1)=d$ (we have
$\hat\alpha_1=\tilde\alpha_1$ if $\tilde\alpha_1\not=\alpha_0$,
$\hat\alpha_1=\alpha_1$ if $\tilde\alpha_1=\alpha_0$), one has $$
\IM \hat\xi_{\hat\alpha_{1}}^{1}=\IM
\xi_{\tilde\alpha_{1}}^{1}<0\; . $$

Conversely, given $(\cA ,\pi_0,\pi_1)$ and $(\zeta_\alpha =
\lambda_\alpha+i\tau_\alpha )_{\alpha\in\cA}$  as above, assume
that $$\IM \xi_{\alpha_{0}}^{0}=\IM \xi_{\alpha_{1}}^{1}\not= 0\;
,  $$ and define $\varepsilon =0$ if $\IM \xi_{\alpha_{1}}^{1}<0$,
$\varepsilon =1$ if $\IM \xi_{\alpha_{0}}^{0}>0$. Set $$ \eqalign{
\hat\zeta_\alpha &= \zeta_\alpha\; , \;\;
\alpha\not=\alpha_\varepsilon\cr \hat\zeta_{\alpha_\varepsilon}
&=\zeta_{\alpha_\varepsilon}+\zeta_{\alpha_{1-\varepsilon}}\; ,
\cr}$$ and define appropriately new combinatorial data; this
operation is the inverse of the one above. Thus the dynamics of
the continued fraction algorithm at the level of suspension is
invertible (on a full measure set) and can be viewed as the
natural extension of the dynamics at the level of i.e.m..

It is clear that the Riemann surfaces $M_\zeta$, $M_{\hat\zeta}$
are canonically isomorphic, and the isomorphism respects the
holomorphic $1$--form and the vertical vector field.

We can also extend the definition of the Zorich algorithm at the
level of suspension data. These accelerated dynamics can actually
be thought of as a first return map of the previous dynamics.
Indeed, in the polyhedral cone of admissible length and suspension
data, consider the polyhedral subcones defined by $$\eqalign{
\cZ_0 &=\{ \lambda_{\alpha_{0}}>\lambda_{\alpha_{1}}\, , \, \IM
\xi_{\alpha_{0}}^{0}>0\}\;, \cr \cZ_1 &=\{
\lambda_{\alpha_{1}}>\lambda_{\alpha_{0}}\, , \, \IM
\xi_{\alpha_{1}}^{1}<0\}\;. \cr}$$ The accelerated dynamics are
the first return map to $\cZ=\cZ_0\sqcup\cZ_1$: this is clear from
the description of the basic step above.

\vskip .5 truecm \noindent {\bf 4.2 The Teichm\"uller flow}

\vskip .5 truecm \noindent Fix combinatorial data $(\cA
,\pi_0,\pi_1)$. Given length data $(\lambda_\alpha
)_{\alpha\in\cA}$ and suspension data $(\tau_\alpha
)_{\alpha\in\cA}$, one defines for $t\in\R$ $$U^t(\lambda ,\tau
)=(e^{t/2}\lambda ,e^{-t/2}\tau)\;  .$$ This flow is called the
{\it Teichm\"uller flow}. Observe that the conditions on the
length data $\lambda_\alpha >0$ and on the suspension data (see
3.1) are preserved under the flow.

It is also obvious that the flow commutes with the basic operation
of the continued fraction algorithm. In particular, the inequality
$\lambda_{\alpha_\varepsilon}>\lambda_{\alpha_{1-\varepsilon}}$ is
preserved.

The surface $M_\zeta$ is canonically equipped with an area form
(coming from $\C$) for which its area is $$A :=\hbox{area}\,
(M_\zeta )=\sum_{\alpha\in\cA}\lambda_\alpha h_\alpha\; . $$ The
area is preserved by the Teichm\"uller flow, and also by the basic
operation of the continued fraction algorithm.

The Lebesgue measure $d\lambda d\tau$ on the domain
$\R^\cA\times\R^\cA$ defined by the restrictions on length and
suspension data is preserved by the Teichm\"uller flow, and by the
basic operation of the continued fraction algorithm.

One now combines the continued fraction algorithm (in Zorich form)
with the Teichm\"uller flow in order to get a version which is
normalized w.r.t.\ scales.

One could decide to normalize by keeping the total length
$\lambda^* =\sum_{\alpha\in\cA}\lambda_\alpha$ constant; actually,
we prefer in the sequel a slightly different normalization, which
leads to simpler formulas.

As in 1.2.1, for
$\lambda_{\alpha_\varepsilon}>\lambda_{\alpha_{1-\varepsilon}}$,
we set $$\eqalign{
\hat\lambda_\alpha &= \lambda_\alpha\; , \;\;
\alpha\not=\alpha_\varepsilon\cr \hat\lambda_{\alpha_\varepsilon}
&=\lambda_{\alpha_\varepsilon}-\lambda_{\alpha_{1-\varepsilon}}\;
. \cr}$$

Define now $$\hat
\lambda^{**}=\sum_{\alpha\in\cA}\hat\lambda_\alpha
=\lambda^*-\lambda_{\alpha_{1-\varepsilon}}\; . $$ Let $(\pi_0
,\pi_1 , \lambda ,\tau )$ belong to the domain $\cZ$ of the Zorich
algorithm, and let $(\bar\pi_0 ,\bar\pi_1 , \bar\lambda ,\bar\tau
)$ be the image. Define $$\eqalign{ t&= t(\lambda) =
2(\log\hat\lambda^{**}-\log\hat{\bar\lambda}^{**})\; , \cr \bar
G(\pi_0 ,\pi_1 , \lambda ,\tau ) &= (\bar\pi_0 ,\bar\pi_1
,U^{t(\lambda)}( \bar\lambda ,\bar\tau ))\; . }$$ The map $\bar G$
is called the normalized step for the natural extension of the
accelerated algorithm.

\vskip .5 truecm \noindent {\bf 4.3 The absolutely continuous
invariant measure}

\vskip .5 truecm \noindent

We already observed that the restriction of Lebesgue measure
$d\lambda d\tau$ to the simplicial cone of admissible length and
suspension data is invariant under both the basic step of the
algorithm and the Teichm\"uller flow.

When we further restrict Lebesgue measure to $\cZ$, we obtain a
measure $m_0$ which is still invariant under Teichm\"uller flow
and is now invariant under the accelerated algorithm.

Observe that the function $t$ used in the definition of $\bar G$
is constant along the orbits of the Teichm\"uller flow. It follows
that the measure $m_0$ is also invariant under $\bar G$.

The area function $A= \sum_{\alpha\in\cA}\lambda_\alpha h_\alpha$
(where $h=-\Omega \tau$) is also invariant under $\bar G$; we
introduce $$\cZ^{(1)} = \cZ\cap\{A\le 1\}\; , $$ and denote by
$m_1$ the restriction of $m_0$ to $\cZ^{(1)}$; it is invariant
under the restriction of $\bar G$ to $\cZ^{(1)}$.

We now project back to the level of i.e.m.\ , i.e.\ of length data
alone: we obtain a map $$G(\pi_0 ,\pi_1, \lambda)=(\bar \pi_0
,\bar\pi_1, e^{t(\lambda )/2}\bar\lambda)$$ and a measure $m_2$
image of $m_1$ under the projection which is invariant under $G$.
As $\hat\lambda^{*}$ is still invariant under $G$, we can
restrict, by homogeneity, the measure $m_2$ to
$\{\hat\lambda^{*}=1\}$ to obtain a measure $m$ invariant under
the restriction of $G$. This is the measure that we are interested
in and that we will now describe.

Let $(\pi_0 ,\pi_1, \lambda)$ be fixed; assume for instance that
$\lambda_{\alpha_0}>\lambda_{\alpha_1}$. Consider in $\tau$--space
the polyhedral cone $$\cU_0 = \{\IM\xi_\alpha^0>0\, ,\,
\forall\alpha\in\cA\, ,\, \IM\xi_\alpha^1<0\,
,\,\forall\alpha\not=\alpha_1\}\; .$$ The density $\chi$ of $m_2$
at $(\pi_0 ,\pi_1, \lambda)$ is equal to the volume of
$\cU_0\cap\{A\le 1\}$. Write $\cU_0$, up to a codimension $1$
subset, as a finite union of disjoint simplicial cones $\cU$. For
each $\cU$, choose a unimodular basis $\tau^{(1)},\ldots
,\tau^{(d)}$ of $\R^\cA$ generating $\cU$ and write
$h^{(j)}=-\Omega\tau^{(j)}$. One has $$\chi_{\pi_0 ,\pi_1}(\lambda
)=(d!)^{-1}\sum_\cU\prod_1^d(\sum_{\alpha\in\cA}\lambda_\alpha
h_\alpha^{(j)})^{-1}\; . \eqno(*)$$ If we set $$\eqalign{
\hat\lambda_{\alpha_0} &=\lambda_{\alpha_0} -\lambda_{\alpha_1}\cr
\hat\lambda_{\alpha} &=\lambda_{\alpha} \, , \, \alpha\not=
\alpha_0\cr \hat h_{\alpha_1} &=h_{\alpha_0} + h_{\alpha_1}\cr
\hat h_{\alpha} &=h_{\alpha} \, , \, \alpha\not= \alpha_1\cr}$$ we
have $$\sum_{\alpha\in\cA}\lambda_\alpha
h_\alpha^{(j)}=\sum_{\alpha\in\cA}\hat\lambda_\alpha \hat
h_\alpha^{(j)}\; . $$ Define $$W_j=\{\alpha\in\cA\, , \,  \hat
h_\alpha^{(j)}\not= 0\}\; .$$

The key property is now the following ([V2],[Z1] see also [Y])

\vskip .5 truecm\noindent {\bf Lemma} {\it For any $X\subset\cA$
with $\emptyset\not= X\not= \cA$, we have} $$\hbox{card}\, \{j\,
,\, W_j\cap X=\emptyset\}+\hbox{card}\, X<d\; .$$

\vskip .3 truecm\noindent When
$\lambda_{\alpha_1}>\lambda_{\alpha_0}$, the only difference is
that we have to start with $$\cU_1 = \{\IM\xi_\alpha^0>0\, ,\,
\forall\alpha\not=\alpha_0\, ,\, \IM\xi_\alpha^1<0\,
,\,\forall\alpha\in\cA\}\; .$$ In the formula $(*)$ above for the
density, set $$\chi_\cU (\lambda )=
\prod_1^d(\sum_{\alpha\in\cA}\lambda_\alpha h_\alpha^{(j)})^{-1}\;
.$$ Up to a constant factor, the density of $m$ on the simplex
$$\Delta =\{\lambda\in\R^\cA\, , \, \hat\lambda_\alpha >0\, , \,
\hat\lambda^{*}=1\}$$ is given by $\sum_\cU\chi_\cU$. One has $$
c^{-1}\le \chi_\cU (\lambda )\prod_{j=1}^d
(\sum_{W_j}\hat\lambda_\alpha)\le c\; .\eqno(1)$$ To control the
size of $\chi_\cU$, we decompose $\Delta$ as follows. Set $$\cN =
\{\vec n=(n_\alpha )_{\alpha\in\cA}\in\N^\cA\, ,\, \min_\alpha
n_\alpha =0\}\; . $$ For $\vec n\in\cN$, $\Delta (\vec n)$ is the
set of $\lambda\in\Delta$ such that $$\eqalign{ \hat\lambda_\alpha
&\ge {1\over 2d}\;\;\hbox{if}\, n_\alpha =0\; , \cr {1\over
2d}2^{1-n_\alpha}>\hat\lambda_\alpha &\ge {1\over
2d}2^{-n_\alpha}\;\;\hbox{if}\, n_\alpha >0\; . \cr}$$ We obtain
thus a partition $$\Delta=\sqcup_\cN\Delta (\vec n)\; , $$ with
the estimate $$c^{-1}\le 2^{\Sigma n_\alpha}\hbox{vol}\,\Delta
(\vec n)\le c\; . \eqno(2)$$ For $\lambda\in\Delta (\vec n)$,
estimate (1) above gives $$c^{-1}\le\chi_\cU (\lambda
)2^{-\Sigma_j\min_{W_j}n_\alpha}\le c\; . \eqno(3)$$

With fixed $\vec n$, let $0=n^0<n^1<\ldots$ be the values taken by
the $n_\alpha$ and $V^i\subset\cA$ the set of indices with
$n_\alpha\ge n^i$. On one side, one has
$$\eqalign{\sum_{\alpha\in\cA}n_\alpha &= \sum_{i\ge
0}n^i(\hbox{card}\, (V^i\setminus V^{i+1}))\cr &= \sum_{i> 0}(n^i
-n^{i-1})\hbox{card}\, V^i\; .\cr}$$ On the other side, let
$\tilde V^i$ be the set of $j$ such that $W_j\subset V^i$; one has
$\min_{W_j}n_\alpha =n^i$ if and only if $j\in\tilde V^i\setminus
\tilde V^{i+1}$ hence $$\eqalign{ \sum_{j=1}^d\min_{W_j}n_\alpha
&= \sum_{i\ge 0}n^i(\hbox{card}\, (\tilde V^i\setminus\tilde
V^{i+1}))\cr &= \sum_{i> 0}(n^i -n^{i-1})\hbox{card}\, \tilde
V^i\; .\cr}$$ By the Lemma above, one has $$\hbox{card}\, \tilde
V^i< \hbox{card}\,  V^i $$ as long as $0<\hbox{card}\, V^i<d$.
This shows that $$\sum_{\alpha\in\cA}n_\alpha
-\sum_{j=1}^d\min_{W_j}n_\alpha \ge |\vec n|_\infty :=
\hbox{Max}_{\alpha\in\cA} n_\alpha\; . $$ The last estimate,
introduced into (2), (3) gives $$ (\hbox{vol}\,\Delta (\vec n))
\hbox{Max}_{\Delta (\vec n)}\chi_\cU \le c 2^{-|\vec n|_\infty}\;
. \eqno (4)$$ The integrability of $\chi_\cU$ over $\Delta$ now
follows from the fact that the number of $\vec n\in\cN$ with
$|\vec n|_\infty=N$ is of order $N^{d-2}$.

If we compare (2) and (4), we obtain $$\hbox{Max}_{\Delta (\vec
n)}\chi_{\pi_0,\pi_1}\le c2^{|\vec n|_1-|\vec n|_\infty}\; ,
\eqno(5)$$ with $|\vec n|_1=\sum_{\alpha\in\cA}n_\alpha$. When
$d=2$, $\chi_{\pi_0,\pi_1}$ is bounded. Assume now $d>2$. From (2)
and (5), one obtains $$
m(\{\chi_{\pi_0,\pi_1}>2^N\})\le\sum_{|\vec n|_1-|\vec
n|_\infty\ge N-c}c2^{-|\vec n|_\infty}\; ; \eqno(6)$$ to have $
|\vec n|_1\ge |\vec n|_\infty+N-c$, one must have $|\vec
n|_\infty\ge{N-c\over d-2}$; an easy computation leads to
$$\eqalignno{m(\{\chi_{\pi_0,\pi_1}>2^N\}) &\le c 2^{-{N\over
d-2}}\; , &(7)\cr \hbox{Leb}\, (\{\chi_{\pi_0,\pi_1}>2^N\}) &\le c
2^{-N{d-1\over d-2}}\; . &(8)\cr}$$

It follows, as $\chi_{\pi_0,\pi_1}$ is bounded from below that we
have, for every Borel set $X$ $$ c^{-1}\hbox{Leb}\, (X)\le m(X)\le
c(\hbox{Leb}\, (X))^{1\over d-1}\; . $$

\vskip .5 truecm \noindent {\bf 4.4 Integrability of $\log\Vert
Z_{(1)}\Vert$}

\vskip .5 truecm \noindent Recall the function $Z_{(1)}$, with
values in $\hbox{SL}\, (\Z^\cA)$, defined in 1.2.4: the sequence
$(\pi_0^{(k)}, \pi_1^{(k)}, \lambda^{(k)})$ given by the Zorich
algorithm satisfies $$\lambda^{(k)} =Z_{(1)}(\pi_0^{(k)},
\pi_1^{(k)}, \lambda^{(k)})\lambda^{(k+1)}\; . $$ Following Zorich
([Z1]) we estimate $\Vert Z_{(1)}\Vert$ w.r.t.\ the absolutely
continuous invariant measure $m$. This will be used in two ways:
\item{$\bullet$} applying Oseledets multiplicative ergodic theorem
in order to prove that conditions (b) and (c) in 1.3 have full
measure;
\item{$\bullet$} as a first step in an induction to prove that
condition (a) in 1.3.1 has full measure.

We use as norm the supremum of the coefficients. For $k\ge 0$,
$\lambda_{\alpha_\varepsilon}>\lambda_{\alpha_{1-\varepsilon}}$,
we have $$\Vert Z_{(1)}\Vert >k\Longleftrightarrow
\hat\lambda_{\alpha_\varepsilon}>k\sum_{\pi_{1-\varepsilon}\alpha
>\pi_{1-\varepsilon}\alpha_\varepsilon}\hat\lambda_{\alpha}\; ; $$
it follows that $$\Vert Z_{(1)}\Vert >(2d)2^{N-1}\Rightarrow
\lambda\in\cup_{|\vec n|_\infty \ge N}\Delta (\vec n )\; , $$
which in turn implies that $$m(\{\Vert Z_{(1)}\Vert > 2^{N}\})\le
cN^{d-2}2^{-N}\; . $$ This is the required estimate; it shows that
$\Vert Z_{(1)}\Vert^\rho$ is $m$--integrable for all $\rho <1$ and
a fortiori that $\log\Vert Z_{(1)}\Vert$ is $m$--integrable.

\vskip .5 truecm \noindent {\bf 4.5 Conditions (b) and (c) have
full measure}

\vskip .5 truecm \noindent As $\log\Vert Z_{(1)}\Vert$ is
$m$--integrable, we can apply Oseledets theorem and obtain the
existence almost everywhere of Lyapunov exponents for the
corresponding cocycle.

The space $\Gamma_s$ is then associated to the negative Lyapunov
exponents. The two estimates in condition (c) are immediate
consequences of the properties of Oseledets decomposition.

For property (b), we recall the result of Veech ([V3]): the
largest Lyapunov exponent is almost everywhere simple. The
existence of a spectral gap follows.

In the end of the section, we will prove that property (a) in
1.3.1 has full measure.

\vskip .5 truecm \noindent {\bf 4.6 The main step}

\vskip .5 truecm \noindent Let $(\cA ,\pi_0,\pi_1)$ be
combinatorial data, $\cD$ the associated Rauzy diagram. For an
i.e.m.\ $T$ satisfying Keane's condition with these data, the
Rauzy--Veech algorithm defines an infinite path
$(\gamma^{(n)}(T))_{n>0}$ in $\cD$, starting at $(\pi_0,\pi_1)$.

Conversely, if $\gamma=(\gamma^{(n)})_{0<n\le N}$ is a finite path
in $\cD$ starting at $(\pi_0,\pi_1)$, we denote by $\Delta (\gamma
)$ the simplex of normalized $T$ in $\Delta (\pi_0,\pi_1)$ such
that $\gamma^{(n)} (T)=\gamma^{(n)}$ for $0<n\le N$. We use here
the old normalization $\{\lambda^*=1\}$.

To such a path $\gamma$ is associated a matrix $Q(\gamma
)\in\hbox{SL}\, (\Z^\cA )$: $$ Q(\gamma ) = V(\gamma^{(1)})\cdots
V(\gamma^{(n)})\; , $$ and we write as before $$Q_\beta (\gamma )
=\sum_{\alpha\in\cA}Q_{\alpha\beta} (\gamma )\; . $$ We have
$$1=\sum_{\alpha\in\cA}\lambda_\alpha^{(0)}=\sum_{\beta\in\cA}
Q_\beta (\gamma )\lambda_\beta$$ (where $\lambda^{(0)}=Q(\gamma
)\lambda^{(N)}$), and it follows that $$\hbox{vol}_{d-1}(\Delta
(\gamma))=[\prod_{\beta\in\cA}Q_\beta (\gamma
)]^{-1}\hbox{vol}_{d-1}(\Delta (\pi_0,\pi_1))\; . $$ Denote by
$(\pi_0^{(N)},\pi_1^{N)})$ the endpoint of $\gamma$, by
$\alpha_0^{(N)}, \alpha_1^{(N)}$ the indices such that
$\pi_\varepsilon^{(N)} (\alpha_\varepsilon^{(N)})=d$. They are the
names of the two arrows going out of $(\pi_0^{(N)},\pi_1^{N)})$.
The conditional probability, for an i.e.m.\ $T$ in $\Delta
(\gamma)$, that the name of $\gamma^{(N+1)}(T)$ is
$\alpha_\varepsilon^{(N)}$ is equal to
$Q_{\alpha_{1-\varepsilon}^{(N)}}(Q_{\alpha_0^{(N)}}+Q_{
\alpha_1^{(N)}})^{-1}$.

Let $1\le D<d$. A segment $(\gamma^{(n)}(T))_{k\le n<l}$ is called
a {\it $D$--segment} if the arrows of the segment take no more
than $D$ distinct names. It is called {\it maximal} if
$(\gamma^{(n)}(T))_{k\le n\le l}$ is not a $D$--segment.

The following proposition is the main step in proving that
condition (a) has full measure.

\vskip .5 truecm \noindent {\bf Proposition}{\it There exist an
integer $l=l(d)$ and a constant $\eta =\eta (d)>0$ with the
following properties. Let $\gamma=(\gamma^{(n)})_{0<n\le N}$ be a
finite path in $\cD$ such that the set $\cA '$ of names of arrows
of $\gamma$ is distinct from $\cA$. Assume that $D=\hbox{card}\,
\cA '>1$. There is a subset $\Delta '(\gamma )$ of $\Delta (\gamma
)$ with $$\hbox{vol}_{d-1}(\Delta '(\gamma
))\ge\eta\hbox{vol}_{d-1}(\Delta (\gamma ))$$ such that, for every
$T\in \Delta '(\gamma )$, there exists $M>N$ with
\item{$\bullet$} the name of $\gamma^{(M)}(T)$ does not belong to
$\cA '$;
\item{$\bullet$} no more than $l$ $(D-1)$--segments are needed to
cover $(\gamma^{(n)})_{N\le n <M}$.}

\vskip .5 truecm \noindent We will first explain how the full
measure estimate for condition (a) follows from the proposition,
and then prove the proposition.

\vskip .5 truecm \noindent {\bf 4.7 Condition (a) has full
measure}

\vskip .5 truecm \noindent For $T\in \Delta (\pi_0,\pi_1)$,
satisfying Keane's condition, and $1\le D<d$, denote by
$Z_{(D)}(T)$ the matrix in $\hbox{SL}\, (\Z^\cA )$ associated to
the initial maximal $D$--segment in $(\gamma^{(n)}(T))_{n>0}$.
Denote by $M_{(D)}(T)$ (resp.\ $M_{(D)}^1(T)$) the minimal number
of $(D-1)$--segments (resp.\ $1$--segments) needed to cover this
initial maximal $D$--segment.

\vskip .5 truecm \noindent {\bf Corollary}{\it Let $N>0$. Except
on a set of measure $\le c 2^{-cN^{1/D}}$, one has} $$\eqalign{
\Vert Z_{(D)}(T)\Vert &\le 2^N\; , \cr M_{(D)}(T) &\le N^{1/D}\; ,
\cr M_{(D)}^1(T) &\le N^{D-1\over D}\; . \cr}$$

\vskip .5 truecm \noindent {\bf Remark}{\it The measure referred
to can be either Lebesgue or the invariant measure $m$: in view of
the last formula of 4.3, it changes only the values of the
constants.}

\vskip .5 truecm \noindent\proof The estimate for $Z_{(1)}$ has
been shown in 4.5. Let us show the estimate for $M_{(D)}$.

Let $\gamma =(\gamma^{(n)})_{0<n\le N}$ be any finite path such
that $(\gamma^{(n)})_{0<n<N}$ is a $(D-1)$--segment but $\gamma$
is not. Apply a first time the proposition in each $\Delta (\gamma
)$. One obtains that $$\hbox{Leb}\, (\{M_{(D)}>l+1\})<1-\eta\; .
$$ We next subdivide the set $\{M_{(D)}>l+1\}$ into simplices
$\Delta (\gamma_1 )$, where $\gamma_1=(\gamma_1^{(n)})_{0<n\le
N_1}$ is a $D$--segment and $(\gamma_1^{(n)})_{0<n< N_1}$ is the
concatenation of $(l+1)$ maximal $(D-1)$--segments. Applying once
again the proposition in each $\Delta (\gamma_1 )$ gives
$$\hbox{Leb}\, (\{M_{(D)}>2l+1\})<(1-\eta)^2\; . $$ Iterating this
process leads to the required estimate for $M_{(D)}$.

We next show by induction on $D$ that $$m(M_{(D)}^1>N^{D-1\over
D})\le c2^{-cN^{1/D}}\; . $$

For $D=2$, one has $M_{(D)}=M_{(D)}^1$; the comparison between $m$
and the Lebesgue measure gives the estimate. Assume $D>2$ and
write $$Z_{(D)}(T)=Z_{(D-1)}(T_0)Z_{(D-1)}(T_1)\cdots
Z_{(D-1)}^*(T_{M-1})$$ with $T_0=T$, $M=M_{(D)}(T)$, and $T_i$ is
obtained from $T_0$ by $n_i$ iterations of the Zorich algorithm
(we have $0=n_0<n_1<n_2<\ldots$); $Z_{(D-1)}^*(T_{M-1})$ denotes
some initial part in the product giving $Z_{(D-1)}(T_{M-1})$.

Neglecting a set of measure $\le c2^{-cN^{1/D}}$, we can assume
that $M\le N^{1/D}$.

By the induction hypothesis, applied with $N'=N^{D-1\over D}$, we
have $$m(M^1_{(D-1)}(T_0)>N^{D-2\over D})\le c 2^{-cN^{1/D}}\; .$$
As the measure $m$ is invariant under the Zorich algorithm, the
same estimate holds when we put instead of $T_0$ any given iterate
$T^{(k)}$ of $T_0$ under the algorithm. Thus we have $$\eqalign{
m(\hbox{Max}_{0\le k<N} M^1_{(D-1)}(T^{(k)})>N^{D-2\over D}) &\le
cN2^{-cN^{1/D}}\cr &\le c'2^{-c'N^{1/D}}\; . \cr}$$ On the other
side, when $$\hbox{Max}_{0\le k<N} M^1_{(D-1)}(T^{(k)})\le
N^{D-2\over D}\; , $$ we have $$n_i\le iN^{D-2\over D}$$ for
$$0\le i\le M-1<N^{1/D}$$ and $$M^1_{(D)}(T)\le N^{D-1\over D}\; .
$$ This proves the estimate for $M^1_{(D)}$.

The estimate on $Z_{(D)}$ is again proven by induction on $D$, the
case $D=1$ having been done in 4.4. Neglecting a set of measure
$c2^{-cN^{1/D}}$, we may assume $M_{(D)}(T)\le N^{1/D}$ and $M^1_{(D)}(T)\le N^{D-1\over
D} $. Write $Z_{(D)}(T)$ as above.

If $\Vert Z_{(D)}(T)\Vert>2^N$, one can find $i\in\{0,1,\ldots
,M-1\}$ such that $$\Vert Z_{(D-1)}(T_i)\Vert >2^{N/M}\ge
2^{N^{D-1\over D}}\; . $$ By the induction hypothesis, we have
$$m\left(\{\Vert Z_{(D-1)}(T_0)\Vert >2^{N^{D-1\over D}}\right)\le
c2^{-cN^{1/D}}$$ and the same estimate holds if we replace $T_0$
by any given $T^{(k)}$. It is sufficient to consider $k\le
N^{{D-1\over D}}$. Again, one has $$cN^{{D-1\over
D}}2^{-cN^{1/D}}\le c'2^{-c'N^{1/D}}\; , $$ and this concludes the
proof of the corollary. \qed

\vskip .5 truecm\noindent The proof that condition (a) has full
measure follows now from a usual Borel--Cantelli argument. Take
$D=d-1$ and write $N=(\kappa \log k)^{d-1}$ with fixed large
$\kappa >0$ and an integer $k\ge 0$. One has $$m\left(\{\Vert
Z_{(d-1)}(T)\Vert >2^{N}\}\right)\le c k^{-c\kappa}\; . $$ If
$\kappa$ is large enough, the right hand term form a converging
series. As $m$ is invariant under the Zorich algorithm, we
conclude that almost surely, the iterates $T^{(k)}$ of $T$ under
the Zorich algorithm satisfy $$\Vert Z_{(d-1)}(T^{(k)})\Vert \le
2^{\kappa (\log k)^{d-1}}\; , $$ for all large $k$.

On the other hand, the exponential rate of growth of the $Q(k)$
(in the Zorich algorithm) is given by the largest Lyapunov
exponent of the Teichm\"uller flow, which is positive.

We conclude that there exists $\kappa_1$ such that almost all
i.e.m.\ $T$ satisfy $$\log\Vert Z_{(d-1)}(k)\Vert\le
\kappa_1[\log\log\Vert Q(k)\Vert]^{d-1}$$ for all large enough
$k$.

\vskip .5 truecm\noindent {\bf Question.}{\it Does one have almost
surely $$\Vert Z_{(d-1)}(k)\Vert=\hbox{O}\, ([\log\Vert
Q(k)\Vert]^C)$$ for some $C>0$?}

\vskip .5 truecm\noindent {\bf 4.8 Proof of the Proposition}

\vskip .5 truecm\noindent Let $\gamma$, $\cA'$, $D$ be as in the
proposition. Let $T$ be an i.e.m.\ in $\Delta (\gamma )$
satisfying Keane's condition. Define, for $n\ge 0$ $$\eqalign{
Q'(n,T) &= \sum_{\alpha\in\cA'} Q_\alpha (n,T) \; , \cr Q_{ext}
(n,T) &= \sum_{\alpha\in\cA\setminus\cA'} Q_\alpha (n,T) \; ,
\cr}$$ where $Q_\alpha (n,T)$ is the shorthand for $Q_\alpha
((\gamma_j(T))_{0\le j\le n})$ (see the beginning of 4.6).

\vskip .5 truecm\noindent {\bf Lemma 1.}{\it If the names of the
arrows $\gamma^{(m)}(T)$ belong to $\cA'$ for $m\le n$, we have
}$$ Q_{ext} (n,T)\le (d-2) Q'(n,T)\; . $$

\vskip .5 truecm\noindent (Recall that, as $1<D<d$, we have $d\ge
3$).

\noindent\proof We start with $Q_\alpha (0,T)=1$ for all
$\alpha\in\cA$. Divide the segment $[1,n]$ into maximal
$1$--segments into which the name of the arrows is the same; let
$[n_i,n_{i+1})$ be such a segment, with arrows of name
$\alpha_i\in\cA'$. The secondary names of these arrows appear with
some periodicity $d_i<d$; moreover, if $n_i>1$, the secondary name
of $\gamma^{(n_i)}$ is $\alpha_{i-1}\in\cA'$; if $n_i=1$, $i=0$,
$n_1>d_0$, the secondary name of $\gamma^{(m)}$ is $\alpha_1$ for
each $m=n_1-kd_0$, $k>d_0$. For $n_i\le m<n_{i+1}$ we have
$$\eqalign{Q_{ext} (m,T) &= Q_{ext} (m-1,T)\; , \cr
Q'(m,T) &= Q'(m-1,T)+Q_{\alpha_i}(n_i-1,T)\; , \cr}$$ if the
secondary name of $\gamma^{(m)}$ is in $\cA'$ and $$\eqalign{
Q_{ext} (m,T) &= Q_{ext} (m-1,T)+Q_{\alpha_i}(n_i-1,T)\; , \cr
Q'(m,T) &= Q'(m-1,T)\; , \cr}$$ otherwise. In each segment except
perhaps the first one, the number of secondary names in
$\cA\setminus\cA'$ does not exceed $(d-2)$ times the number of
secondary names in $\cA'$. In the first segment, we write
$n_1=kd_0+n_1'$, $0<n_1'\le d_0$; again the number of secondary
names in $\cA\setminus\cA'$ does not exceed $(d-2)$ times the number of
secondary names in $\cA'$ in the subsegment $[n_1',n_1)$. Finally
we have for $0\le m<n_1'$ that $Q'(m,T)\ge D\ge 2$, $Q_{ext}
(m,T)\le Q_{ext} (0,T)+m\le d-D+d_0-1\le 2d-4$ and the estimate of
the lemma follows. \qed

\vskip .5 truecm\noindent Let $1\le D_1\le D$, $n\ge 0$, $C_1>0$.
We say that $T\in \Delta (\gamma )$ is {\it
$(D_1,n,C_1)$--balanced} if we have $$Q_\alpha (n,T)\ge C_1^{-1}
Q'(n,T)$$ for at least $D_1$ indices $\alpha\in\cA'$. The property
only depends on the path $(\gamma^{(m)}(T))_{0<m\le n}$ and we
will also say that this path is {\it $(D_1,n,C_1)$--balanced}.
Clearly, any $T$ is $(1,n,D)$--balanced (for all $n\ge 0$).

\vskip .5 truecm\noindent {\bf Lemma 2.}{\it Assume that $\gamma$
is $(D,n,C_0)$--balanced, for some constant $C_0>0$. Then we can find
$\Delta '(\gamma )\subset \Delta (\gamma )$ satisfying the conclusions
of the proposition, with $l=l(d)$ and $\eta =\eta (d,C_0)$.}

\vskip .5 truecm\noindent\proof Let $\gamma '=(\gamma^{(n)})_{0<
n\le M}$ be an extension of $\gamma$ with minimal length such that
the name $\alpha$ of $\gamma^{(M)}$ is not in $\cA'$. Then $M-N$
is bounded by the diameter of $\cD$, i.e.\ in terms of $d$ only.
Therefore there exists $C_*=C_*(d)$ such that $\gamma '$ is
$(D,M-1,C_*C_0)$--balanced; moreover, the path $\gamma''=
(\gamma^{(n)})_{0< n< M}$ satisfies $$ \hbox{Vol}_{d-1}\, \Delta
(\gamma'')\ge\eta'' \hbox{Vol}_{d-1}\, \Delta (\gamma )\; , $$
with $\eta''=\eta''(C_0,d)$. Then, for all $\beta\in\cA'$, we have
$$\eqalign{ Q_\alpha (M-1) &\le Q_{ext}(M-1)\cr &\le
(d-2)Q'(M-1)\cr &\le (d-2)C_*C_0Q_\beta (M-1)\; , \cr}$$ and
therefore $$\hbox{Vol}_{d-1}\, \Delta
(\gamma')\ge\eta'\hbox{Vol}_{d-1}\, \Delta (\gamma'')\; , $$ with
$\eta'=(1+(d-2)C_*C_0)^{-1}$. We take $\eta =\eta'\eta''$, $\Delta
'(\gamma )=\Delta (\gamma')$. Finally $l$ is bounded because $M-N$
is bounded. \qed

\vskip .5 truecm\noindent When $\gamma$ is only $(\tilde
D,N,\tilde C)$--balanced for some $\tilde D<D$, the strategy will
be to extend $\gamma$ without losing volume in order to obtain a
more balanced path; at the end we should be able to apply Lemma 2
(unless we have already found $\Delta '(\gamma )$).

We therefore assume that $\gamma$ is $(\tilde
D,N,\tilde C)$--balanced. This is certainly satisfied with $\tilde
D=1$, $\tilde C=D$. Denote by $\tilde \cA$ the set of $\alpha$
such that $$Q_\alpha (N)\ge\tilde C^{-1}Q'(N)\; .$$ The first step
is to extend $\gamma$ to a path $\gamma'=(\gamma^{(n)}(T))_{0<n\le N'}$
of minimal length such that the name of $\gamma^{(N')}$ is not in
$\tilde\cA$. When $N'=N+1$, there might be two choices for
$\gamma^{(N')}$and we choose the one which gives the largest
volume to $\Delta (\gamma ')$.

In any case, an argument completely similar to the one in the
proof of Lemma 2 leads to the estimate $$\hbox{Vol}_{d-1}\, \Delta
(\gamma ')\ge\eta' \hbox{Vol}_{d-1}\, \Delta (\gamma )\; , $$ with
a constant $\eta'=\eta'(d,\tilde C)$.

If the name of $\gamma^{(N')}$ does not belong to $\cA'$, we can
take as in Lemma 2 $\Delta'(\gamma )=\Delta (\gamma ')$ and the
proof of the proposition is over. We now assume that the name of
$\gamma^{(N')}$ belongs to $\cA'\setminus\tilde\cA$.

The subset $\Delta'(\gamma )$ of $\Delta(\gamma )$ we are looking
for will be contained in $\Delta (\gamma ')$. Observe that there
exists $C_*=C_*(d)$ such that $\gamma'$ is $(\tilde D,N',C_*\tilde
C)$--balanced.

\vskip .3 truecm\noindent {\bf Case A:} In the loop of arrows of
the same name which starts with $\gamma^{(N')}$, no secondary name
belongs to $\cA'\setminus\tilde\cA$.

Let $\alpha$ be the name of $\gamma^{(N')}$, $\beta_0,\ldots
,\beta_{r-1}$ being the successive secondary names in the loop.
Let $k>0$, that we write $k=rl+m$, $0\le m<r$. Let $\gamma_1(k)$
be the path extending $\gamma'$ such that
\item{$\bullet$} the name of $\gamma_1(k)^{(n)}$ is $\alpha$ for
$N'<n<N'+k:=N_1(k)$;
\item{$\bullet$} the name of $\gamma_1(k)^{(N_1(k))}$ is
$\beta_m$.

Observe that it follows immediately from the definition of
$R_0,R_1$ in 1.2.1 that the indices $\beta_0,\ldots
,\beta_{r-1}$ are distinct. Therefore, we will have, for $0\le
k_1=rl_1+m_1<k$:
$$\eqalign{ Q_{\beta_j}(N'+k_1) &=
Q_{\beta_j}(N'-1)+\cases{l_1Q_\alpha (N'-1) &if $m_1<j$, \cr
(l_1+1)Q_\alpha (N'-1) &if $m_1\ge j$, \cr}\cr
Q_\alpha (N'+k_1) &= Q_\alpha (N'-1)\; , \cr}$$ and also
$$\eqalign{Q_{\beta_j}(N'+k) &=
Q_{\beta_j}(N'+k-1)\; , \cr
Q_\alpha (N'+k) &= Q_\alpha (N'-1)+Q_{\beta_m}(N'+k-1)\; . \cr}$$

For any $k>0$, the extension from $\gamma$ to $\gamma_1(k)$ is
covered by the same number of $(D-1)$ segments, which is bounded
in terms of $d$ only.

For those $k$ such that $\beta_m\notin\cA'$, we include $\Delta
(\gamma_1(k))$ in $\Delta'(\gamma )$.

The formulas for the volumes give $$\hbox{Vol}_{d-1}\left(\left[
\cup_{0<k_1<k}\Delta
(\gamma_1(k_1))\right]^c\right)={\prod Q_{\beta_j}(N')\over\prod
Q_{\beta_j}(N'+k-1)}\hbox{Vol}_{d-1}\Delta (\gamma' )\; . $$

We keep for further consideration all $\gamma_1(k)$ with
$$kQ_\alpha (N'-1)\le Q'(N'-1)\; .$$

The formula above shows that together they will fill a definite
proportion of $\Delta (\gamma' )$.

We also see that when $\beta_m\in\tilde\cA$, $\gamma_1(k)$ will be
$(\tilde D+1,N_1(k),C_1)$--balanced, with $C_1$ depending only on
$d$. For each such $\gamma_1(k)$, we either apply Lemma 2 (if
$\tilde D+1=D$) or repeat the discussion, with $\gamma_1(k)$ in
the place of $\gamma$, from a better starting hypothesis.

\vskip .3 truecm\noindent {\bf Case B:} The complement of case A.

For an i.e.m.\ in $\Delta (\gamma' )$ satisfying Keane's
condition, we consider the three mutually exclusive possibilities:

\item{$\bullet$} $T$ {\it is of type I} if there exists $N_1\ge
N'$ such that all arrows $\gamma^{(n)}(T)$, $N'\le n\le N_1$, have
names in $\cA'\setminus\tilde\cA$, and we have $$\sum_{\alpha\in
\cA'\setminus\tilde\cA}Q_\alpha (N_1,T)\ge Q'(N',T)\; . $$
We take a minimal such $N_1$.
\item{$\bullet$} $T$ {\it is of type II} (respectively {\it of
type III}) if it is not of type I and the first name of an arrow
$\gamma^{(n)}(T)$, $n>N'$, which does not belong to $\cA'\setminus\tilde\cA$
belongs to $\tilde\cA$ (resp.\ to $\cA\setminus\cA'$).

We deal separately with the three types.

\vskip .3 truecm\noindent {\bf a)} All $T$ of type III will be
contained in $\Delta'(\gamma )$; for such a $T$, $M$ is the first
integer $>N'$ for which the name does not belong to
$\cA'\setminus\tilde\cA$. Observe that the segment
$(\gamma^{(n)}(T))_{N'\le n<M}$ is a $(D-1)$--segment because
$\hbox{card}\, (\cA'\setminus\tilde\cA )<D$. As $N'-N$ is bounded
in function of $d$ only, the number of $(D-1)$--segments needed to
cover $(\gamma^{(n)}(T))_{N\le n<M}$ is bounded in terms of $d$
only.

\vskip .3 truecm\noindent {\bf b)} Assume that $T$ is of type II.
Let $N_1$ be the smallest integer $n>N'$ such that the name of
$\gamma^{(n)}(T)$ does not belong to $\cA'\setminus\tilde\cA$;
this name belongs to $\tilde\cA$. Let
$\gamma_1=(\gamma^{(n)}(T))_{0<n\le N_1}$. When $T$ varies among
i.e.m.\ 's of type II, the $\gamma_1$ form an at most countable
collection such that the corresponding simplices $\Delta (\gamma_1
)$ have disjoint interiors (and are contained in $\Delta (\gamma'
)$). Every $T_1$ belonging to some $\Delta (\gamma_1
)$ is also of type II. We claim that every $\gamma_1$ is
$(D_1,N_1,C_1)$--balanced with $D_1>\tilde D$ and $C_1=C_1(\tilde
C,d)$ (see the proof below). As for type III, the number of
$(D-1)$--segments needed to cover $\gamma^{(n)}_1$, $N\le n<N_1$,
is bounded in terms of $d$ only.

\vskip .3 truecm\noindent {\bf c)} Assume that $T$ is of type I.
With $N_1$ minimal as in the definition of type I, take
$\gamma_1=(\gamma^{(n)}(T))_{0< n\le N_1}$. When $T$ varies among
i.e.m.\ 's of type I, the $\gamma_1$ form again an at most
countable collection for which the corresponding simplices $\Delta
(\gamma_1 )$ have disjoint interiors (and are contained in $\Delta
(\gamma' )$). Every $T_1$ belonging to some $\Delta (\gamma_1)$ is
also of type I. We claim that every $\gamma_1$ is
$(D_1,N_1,C_1)$--balanced with $D_1>\tilde D$ and $C_1=C_1(\tilde
C,d)$ (see the proof below). The number of $(D-1)$--segments needed
to cover $\gamma^{(n)}_1$, $N\le n<N_1$,
is bounded in terms of $d$ only.

The discussion above leads in case B to a countable partition (up
to a codimension one subset) of $\Delta (\gamma' )$ into
subsimplices of type III which will be included in $\Delta'
(\gamma )$ and simplices $\Delta (\gamma_1)$ (of type I or II)
which satisfy the same hypotheses than $\Delta (\gamma )$ but are
better balanced (i.e.\ $D_1>\tilde D$); when $D_1=D$, we can apply
Lemma 2 to $\gamma_1$; when $D_1<D$, we repeat the discussion with
$\gamma_1$ instead of $\gamma$. The process stops in less than $D$
steps and gives the conclusion of the proposition. \qed

\vskip .5 truecm\noindent {\it Proof of the claim for type II.} As
$T$ is not of type I, we have
$$\sum_{\alpha\in\cA'\setminus\tilde\cA} Q_\alpha
(N_1-1,T)<Q'(N',T)\; . $$

Let us consider a maximal $1$--segment contained in
$(\gamma^{(n)}(T))_{N'\le n< N_1}$. As we are not in case A, there
is a definite proportion, depending only of $d$, of secondary
names which belong to $\cA'\setminus\tilde\cA$.  This implies that
we must have $$Q'(N_1-1,T)\le C_1'Q'(N',T)\; , $$ with $C_1'$
depending only on $d$. On the other hand, if $\alpha\in
\cA'\setminus\tilde\cA$ and $\beta\in\tilde\cA$ are the names of
$\gamma^{(N_1-1)}(T)$, $\gamma^{(N_1)}(T)$ respectively, we have
$$\eqalign{ Q_\alpha (N_1) &= Q_\alpha (N_1-1)+Q_\beta (N_1-1)\cr
&\ge (C_*\tilde C)^{-1}Q'(N',T)\; . \cr}$$ It follows that
$\gamma_1$ is $(\tilde D+1,N_1,C_1)$--balanced with
$C_1=C_1'C_*\tilde C$. \qed

\vskip .5 truecm\noindent {\it Proof of the claim for type I.} By
definition of $N_1$, we have again
$$\sum_{\alpha\in\cA'\setminus\tilde\cA} Q_\alpha (N_1-1,T)<
Q'(N',T)\; , $$ and it follows again that $$Q'(N_1-1,T)\le
C_1'Q'(N',T)\; . $$ By definition of $N_1$, we have now $$
\sum_{\alpha\in\cA'\setminus\tilde\cA} Q_\alpha (N_1,T)\ge
Q'(N',T)\; , $$ and it follows that $\gamma_1$ is $(\tilde
D+1,N_1,C_1)$--balanced with $C_1$ depending only on $\tilde C$
and $d$. \qed

The proof of the proposition, and therefore also of the full
measure statement, is now complete.

\vfill\eject
\noindent {\bf Appendix A. Roth--type conditions in a concrete
family of interval
exchange maps }

\vskip .5 truecm \noindent {\bf A.1} Let $\cA =(A,B,C,D)$. The
Rauzy diagram of the pair $(\pi_0,\pi_1)=\pmatrix{
A & B & C & D \cr D & C & B & A\cr}$ is indicated in 1.2.2. The
suspension of an i.e.m.\ with these combinatorial data leads to an
holomorphic $1$--form with a double zero on a genus two surface.

In this diagram, we define for $n\ge 0$ a loop $\gamma (n)$ based
at $(\pi_0,\pi_1)$ by asking that the names of the successive
arrows should be $D^2CDA^2B^nA$. The product of the $V$ matrices
around this loop is $$M(n)=\pmatrix{1 & 1 & 1 & 1\cr n & n+1 & 0 &
0 \cr 0 & 0 & 2 & 1\cr n+1 & n+2 & 2 & 2\cr}$$ with characteristic
polynomial $$\chi_n(X)=X^4-(n+6)X^3+(3n+10)X^2-(n+6)X+1$$ Setting
$U=X+X^{-1}$ leads to $$\chi_n(X)=X^2(U^2-(n+6)U+3n+8)\; . $$ The
eigenvalues of $M(n)$ are thus given by
$$\lambda+\lambda^{-1}=U^\pm :={1\over 2}(n+6\pm\sqrt{n^2+4})\; .
$$ The case $n=0$ is degenerate, with $U^+=4$, $U^-=2$. When
$n>0$, both $U^+$, $U^-$ are $>2$; we will denote the eigenvalues
by $\lambda^+_u>\lambda^-_u (>1)>\lambda^-_s>\lambda^+_s$, by
$e_u^+,e_u^-,e_s^-,e_s^+$ the corresponding eigenvectors of the
{\it transposed} matrix $^tM(n)$.

The eigenvector associated to the eigenvalue $\lambda$ is
proportional to $$((\lambda -1)(\lambda^2-4\lambda +2),
\lambda^3-4\lambda^2+3\lambda-1,\lambda (\lambda-1),
(\lambda-1)^2)\; . $$

\vskip .5 truecm \noindent {\bf A.2} As $n\rightarrow +\infty$,
one has $$\eqalign{ \lim U^+-(n+3) &= \lim\lambda_u^+-(n+3)=0\; ,
\cr \lim U^- &=3\; ,\;\;\lim\lambda_u^-=G:={\sqrt{5}+3\over 2}\; .
\cr}$$ One can also choose eigenvectors to obtain: $$\eqalign{
\lim e_u^+ &= E_u^+ := (1,1,0,0)\cr
\lim e_u^- &= E_u^- := (-1,-1,G-1,1)\cr
\lim e_s^- &= E_s^- := (-1,-1,G^{-1}-1,1)\cr
\lim e_s^+ &= E_s^+ := (2,1,0,-1)\; . \cr}$$
These four limit vectors form a basis of $\R^4$ in which we
rewrite $^tM(n)$:
$$\eqalign{ ^tM(n)E^+_u
&=(n+3)E^+_u-E^+_s+{1\over\sqrt{5}}(E^-_u-E^-_s)\; , \cr
^tM(n)E^+_s
&=E^+_u\; , \cr ^tM(n)E^-_u
&=G(E^+_u+E^-_u)\; , \cr ^tM(n)E^-_s
&=G^{-1}(E^+_u+E^-_s)\; . \cr}$$ For the corresponding
coordinates, this gives $$\eqalign{ X^+_u &=
(n+3)x^+_u+x^+_s+Gx_u^-+G^{-1}x_s^-\; , \cr X^+_s &=
-x^+_u\; , \cr X^-_u &=
{1\over\sqrt{5}}x^+_u+Gx_u^-\; , \cr X^-_s &=
-{1\over\sqrt{5}}x^+_u+G^{-1}x_u^-\; . \cr}$$

\vskip .5 truecm \noindent {\bf A.3} The following two lemmas
express that for $n\ge 4$ certain cone conditions are satisfied.

\vskip .5 truecm \noindent {\bf Lemma 1.}{\it For $n\ge 4$, $x^+_u\ge
\hbox{Max}\,(|x^+_s|,|x^-_u|,|x^-_s|)$ one has $$ X^+_u\ge \hbox{Max}\,\left((n-1)|X^+_s|,
(n-1)|X^-_s|,{10-3\sqrt{5}\over 3}|X^-_u|\right)\; , $$ and}
$$X^+_u\ge(n-1)x^+_u\; . $$

\vskip .5 truecm\noindent\proof As $G^{-1}+{1\over\sqrt{5}}<1$ and
$G+G^{-1}=3$, we have $$\hbox{Max}\,\left((n-1)|X^+_s|,
(n-1)|X^-_s|,(n-1)|x^+_u|\right)\le X^+_u\; . $$ If $x^-_u\ge 0$,
one has $$ X^+_u\ge (n+1)x^+_u+Gx^-_u\ge 2X^-_u$$ because
$n+1-2/\sqrt{5}>G$ for $n\ge 4$.

For $x^-_u<0$, one has $X^+_u\ge (n+1)x^+_u-G|x_u^-|$. On one hand
$$(n+1)x^+_u-G|x_u^-|\ge (n+1)\sqrt{5}X_u^-\; , $$ on the other
$$(n+1)x^+_u-G|x_u^-|\ge -\gamma X^-_u$$ as soon as
$n+1+\gamma/\sqrt{5}\ge (\gamma +1)G$, which allows to take
$\gamma={10-3\sqrt{5}\over 3}$ for $n\ge 4$.\qed

\vskip .5 truecm \noindent {\bf Lemma 2.}{\it For $n\ge 4$, $
\hbox{Max}\,(|x^+_u|,|x^-_u|)\ge \hbox{Max}\,(|x^+_s|,|x^-_s|)$, one
has}
$$\eqalign{
\hbox{Max}\,(|X^+_u|,|X^-_u|)&\ge\left(G-{1\over\sqrt{5}}\right)
\hbox{Max}\,(|X^+_s|,|X^-_s|)\; ,\cr
\hbox{Max}\,(|X^+_u|,|X^-_u|)&\ge\left(G-{1\over\sqrt{5}}\right)
\hbox{Max}\,(|x^+_u|,|x^-_u|)\; .\cr}$$

\vskip .5 truecm\noindent\proof When $|x^+_u|\ge |x^-_u|$, this
follows from Lemma 1. If $0\le |x^+_u|\le x_u^-$, one has
$$X_u^-\ge
\left(G-{1\over\sqrt{5}}\right)x^-_u\ge\left(G-{1\over\sqrt{5}}\right)
\hbox{Max}\,(|X^+_s|,|X^-_s|)$$ because
$G^{-1}+{1\over\sqrt{5}}<1$. \qed

\vskip .5 truecm\noindent One should observe in Lemma 1 that ${10-3\sqrt{5}\over
3}>1$ and in Lemma 2 that $G-{1\over\sqrt{5}}>1$.

\vskip .5 truecm \noindent {\bf Lemma 3.}{\it Equip $\R^\cA$ with
the sup norm. Then, for any integers $n_1,\ldots ,n_k>0$, we have}
$$\prod_{i=1}^k(n_i+1)\le\Vert ^tM(n_k)\cdots
^tM(n_1)\Vert_\infty\le\prod_{i=1}^k(2n_i+4)\; . $$

\vskip .5 truecm\noindent\proof The upper bound follows from
$\Vert ^tM(n)\Vert\le 2n+4$ for $n>0$, the lower bound from the
fact that $$^tM(n)(1,1,0,0)-(n+1)(1,1,0,0)$$ is a non negative
vector. \qed

\vskip .5 truecm \noindent {\bf A.4} Let $\Sigma$ be the set of
sequences $(n_i)_{i>0}$ of integers $\ge 4$. To each sequence in
$\Sigma$ we associate the infinite path $\gamma (n_1)\gamma
(n_2)\cdots$ starting at $(\pi_0,\pi_1)$. The cone property of
Lemma 1 guarantees that there is exactly one i.e.m.\ satisfying
Keane's condition associated with this path. On the space $\Gamma$
of functions constant on each $j_0(I_\alpha )$, we have a complete
filtration: the space $\Gamma_s$ has dimension $2$ according to
Lemma 2, contains the line $\R\delta$ (where $\delta$ is the
displacement vector) and is contained in the hyperplane $\Gamma_*$
of zero mean.

Therefore conditions (b) and (c) in 1.3.2, 1.3.3 are automatically
satisfied. Condition (a) is equivalent, in view of Lemma 3, to
$$\log n_k=\hbox{o}\, \left(\sum_{i<k}\log n_i\right)\; . $$

\vfill\eject
\noindent {\bf Appendix B. A non--uniquely ergodic interval
exchange map satisfying condition (a) }

\vskip .5 truecm \noindent {\bf B.1} Let $m,n,p$ be non negative
integers. In the Rauzy diagram of the pair
$(\pi_0,\pi_1)=\pmatrix{
A & B & C & D \cr D & C & B & A\cr}$ (cf.\ 1.2.2), consider the
loop $\gamma_0(m,n,p)$ based at
$(\pi_0,\pi_1)$ such that the names of the successive arrows are
$$ D^{3m+1}BC^nBDC^pD.$$ We also consider the
dual loop $\gamma_1(m,n,p)$ which is deduced from
$\gamma_0(m,n,p)$ by means of the canonical involution and whose
arrows have names $$ A^{3m+1}CB^nCAB^pA.$$ Given three
sequences $(m_k)_{k\ge 0}$, $(n_k)_{k\ge 0}$ and $(p_k)_{k\ge 0}$
we also consider the infinite path $\Gamma$, based at
$\pmatrix{A&B&C&D\cr D&C&B&A\cr}$ which is obtained composing $$
\gamma_0(m_0,n_0,p_0)\gamma_1(m_1,n_1,p_1)\cdots\gamma_0(m_{2k},n_{2k},p_{2k})
\gamma_1(m_{2k+1},n_{2k+1},p_{2k+1})\cdots$$
The matrix $Z_0(m,n,p)$ associated to $\gamma_0(m,n,p)$ is $$
\pmatrix{ 1&0&0&0\cr 0&2&p+2&p+1\cr 0&n&(n+1)(p+1)&p(n+1)\cr
m+1&m(n+2)+1&m(n+2)(p+1)+m+1&pm(n+2)+m+1\cr}$$ where the vectors
of the canonical basis of $\R^\cA$ are ordered alphabetically.
Analogously the matrix $Z_1(m,n,p)$ associated to
$\gamma_1(m,n,p)$ is $$\pmatrix{
pm(n+2)+m+1&m(n+2)(p+1)+m+1&m(n+2)+1&m+1\cr p(n+1)&(n+1)(p+1)
&n&0\cr p+1&p+2&2&0\cr 0&0&0&1\cr}$$ We set
$$Q(k)=Z_0(m_0,n_0,p_0)Z_1(m_1,n_1,p_1)\cdots Z_\varepsilon
(m_{k-1},n_{k-1},p_{k-1})\; ,$$ with
$k-1\equiv\varepsilon\,\hbox{mod}\, 2$. We denote
$e_A(k),e_B(k),e_C(k),e_D(k)$ the column vectors of $Q(k)$.

\vskip .5 truecm \noindent {\bf B.2} Let $m_0=0$ and choose
$n_0\gg 1$. The integer $p_0$ will be chosen later but it will be
such that $p_0\ge n_0$. One has $$\eqalign{ e_A(1)
&=\;^t\pmatrix{1 & 0 & 0 & 1\cr}\, , \cr e_B(1) &=
n_0\left[\;^t\pmatrix{0 & 0 & 1 & 0\cr} + \hbox{O}\,
(n_0^{-1})\right]\, , \cr e_C(1) &= n_0p_0\left[\;^t\pmatrix{0 & 0
& 1 & 0\cr} + \hbox{O}\, (n_0^{-1})\right]\, , \cr e_D(1) &=
n_0p_0\left[\;^t\pmatrix{0 & 0 & 1 & 0\cr} + \hbox{O}\,
(n_0^{-1})\right]\, . \cr}$$ We determine then
$m_1,p_0,n_1,m_2,\ldots ,m_k,p_{k-1},n_k,m_{k+1},p_k,\ldots $
through the following formulas: $$\eqalign{\Pi_0 &:=n_0\, ,\cr
\Pi_1 &:= m_1\Pi_0^{-1}=n_0^2\, ,\cr \Pi_2 &:=
p_0\Pi_1^{-1}=(n_0+1)^2\, , \cr &\vdots \cr \Pi_{3l} &:=
n_l\Pi_{3l-1}^{-1}=(n_0+3l-1)^2\, , \cr \Pi_{3l+1} &:=
m_{l+1}\Pi_{3l}^{-1}=(n_0+3l)^2\, ,\cr \Pi_{3l+2} &:=
p_l\Pi_{3l+1}^{-1} = (n_o+3l+1)^2\, , \cr &\vdots
\cr}\quad\eqalign{m_1 &:= n_0^3\, , \cr p_0 &:= (n_0+1)^2\Pi_1 =
n_0^2(n_0+1)^2\, , \cr n_1 &:= (n_0+2)^2\Pi_2=(n_0+1)^2(n_0+2)^2\,
, \cr &\vdots \cr m_{l+1} &:= (n_0+3l)^2\Pi_{3l}\, , \cr p_l &:=
(n_0+3l+1)^2\Pi_{3l+1}\, , \cr n_{l+1} &:=
(n_0+3l+2)^2\Pi_{3l+2}\, , \cr &\vdots \cr}$$

Thus one has, for $l\ge 0$ $$\eqalign{p_l &=
(n_0+3l)^2(n_0+3l+1)^2\, , \cr n_{l+1} &=
(n_0+3l+1)^2(n_0+3l+2)^2\, ,\cr m_{l+2} &=
(n_0+3l+2)^2(n_0+3l+3)^2\, ,\cr}$$ and also $m_1=n_0^3$. For all
$k\ge -1$ we set $$ c_k=n_0^3\left[{(n_0+k)!\over n_0!}\right]^2\,
$$ so that one has $n_0=c_{-1}$ and $$ \eqalign{ c_1 &= n_0p_0\,
,\cr c_2 &= m_1n_1\, ,\cr c_4 &= m_1n_1p_1\, ,\cr c_5 &=
m_2n_2n_0p_0\, ,\cr c_7 &= m_2n_2p_2n_0p_0\, , \cr c_8
&=m_3n_3m_1n_1p_1\, , \cr &\vdots \cr c_{10} &=m_3n_3p_3m_1n_1p_1\, , \cr
&\vdots \cr}$$

Let us check by induction that, setting $c_{-2}=1$, one has for
$l\ge 0$: $$ \eqalign{ e_A (2l-1) &=c_{6l-8}\left[\;^t\pmatrix{1 &
0 & 0 & 1\cr}+\hbox{O}\, (n_0^{-1})\right]\, , \cr e_B (2l-1)
&=c_{6l-7}\left[\;^t\pmatrix{0 & 0 & 1 & 0\cr}+\hbox{O}\,
(n_0^{-1})\right]\, , \cr e_C (2l-1)
&=c_{6l-5}\left[\;^t\pmatrix{0 & 0 & 1 & 0\cr}+\hbox{O}\,
(n_0^{-1})\right]\, , \cr e_D (2l-1)
&=c_{6l-5}\left[\;^t\pmatrix{0 & 0 & 1 & 0\cr}+\hbox{O}\,
(n_0^{-1})\right]\, , \cr e_D (2l) &=c_{6l-5}\left[\;^t\pmatrix{0
& 0 & 1 & 0\cr}+\hbox{O}\, (n_0^{-1})\right]\, , \cr e_C (2l)
&=c_{6l-4}\left[\;^t\pmatrix{1 & 0 & 0 & 1\cr}+\hbox{O}\,
(n_0^{-1})\right]\, , \cr e_B (2l) &=c_{6l-2}\left[\;^t\pmatrix{1
& 0 & 0 & 1 \cr}+\hbox{O}\, (n_0^{-1})\right]\, , \cr e_A (2l)
&=c_{6l-2}\left[\;^t\pmatrix{1 & 0 & 0 & 1\cr}+\hbox{O}\,
(n_0^{-1})\right]\, , \cr } $$ We have already checked the first
four relations for $l=1$. Assume that the first four relations are
verified for a given value of $l$. Then $$ e_D
(2l)=e_D(2l-1)+(m_{2l-1}+1)e_A(2l-1)\, , $$ with $$
m_{2l-1}c_{6l-8}=c_{6l-6}=c_{6l-5}(n_0+6l-5)^{-2}\, . $$ Moreover
$$ e_C(2l)=2e_C(2l-1)+n_{2l-1}e_B(2l-1)+[m_{2l-1}(n_{2l-1}+2)+1]
e_A(2l-1)\, , $$ with $$ \eqalign{ m_{2l-1}n_{2l-1}c_{6l-8} &=
c_{6l-4}\, , \cr n_{2l-1}c_{6l-7} &= (n_0+6l-6)^{-2}c_{6l-4}\, ,
\cr 2c_{6l-5} &=2(n_0+6l-4)^{-2}c_{6l-4}\, , \cr
(2m_{2l-1}+1)c_{6l-8} &=\hbox{O}\,((n_0+6l)^{-2}) c_{6l-4}\, , \cr
} $$ and $$ \eqalign{ e_B(2l) &=
(p_{2l-1}+2)e_C(2l-1)+(n_{2l-1}+1)(p_{2l-1}+1) e_B(2l-1)\cr
&\phantom{=}+(m_{2l-1}(n_{2l-1}+2)(p_{2l-1}+1)+m_{2l-1}+1)
e_A(2l-1)\cr} $$ with $$ \eqalign{ (p_{2l-1}+2)c_{6l-5}
&=\hbox{O}\,((n_0+6l)^{-2} c_{6l-2})\, , \cr
(n_{2l-1}+1)(p_{2l-1}+1)c_{6l-7} &= \hbox{O}\,((n_0+6l)^{-2}
c_{6l-2})\, , \cr [1+m_{2l-1}(n_{2l-1}+2p_{2l-1}+1)]c_{6l-8}
&=\hbox{O}\,((n_0+6l)^{-4} c_{6l-2})\, , \cr
m_{2l-1}n_{2l-1}p_{2l-1}c_{6l-8} &=c_{6l-2}\, .\cr} $$ The formula
for $e_A(2l)$ is completely similar.

Since one has
$$
\prod_{k\ge 0}[1+\hbox{O}\,((n_0+k)^{-2})]-1=
\hbox{O}\,(n_0^{-1})\, ,
$$
one gets the four last relations. Taking into account the canonical
involution one can analogously obtain the first four relations.

\vskip .5 truecm \noindent {\bf B.3} The decomposition of the
infinite path $\Gamma$ into loops $\gamma_0(m_{2k},n_{2k},p_{2k})$
and $\gamma_1(m_{2k+1},n_{2k+1},p_{2k+1})$ is nothing else than
the decomposition for the accelerated Zorich algorithm.
One has
$$
\Vert Z_\varepsilon (m_k,n_k,p_k)\Vert \sim (n_0+3k)^{12}\, ,
$$
(with $\varepsilon\equiv k\,\hbox{mod}\, 2$), and
$$
\Vert Q(k)\Vert \sim c_{3k-2}=n_0^3\left[{(n_0+3k-2)!\over
n_0!}\right]^2\, .
$$
Thus one obtains
$$
\Vert Z_\varepsilon (m_k,n_k,p_k)\Vert=\hbox{o}\, \left(\left[
\log \Vert Q(k)\Vert\right]^{12}\right)
$$
and the first condition in the definition of Roth type interval exchange map
is (by far) satisfied.

\vskip .5 truecm \noindent {\bf B.4} From the formula and
estimates of Section A.2.2 one gets
$$
\eqalign{
{e_D(2l)\over\Vert e_D(2l)\Vert_1} &= {e_D(2l-1)\over\Vert e_D(2l-1)\Vert_1}
+\hbox{O}\,((n_0+6l)^{-2})\, ,\cr
{e_C(2l)\over\Vert e_C(2l)\Vert_1} &= {e_A(2l-1)\over\Vert e_A(2l-1)\Vert_1}
+\hbox{O}\,((n_0+6l)^{-2})\, ,\cr
{e_B(2l)\over\Vert e_B(2l)\Vert_1} &= {e_A(2l-1)\over\Vert e_A(2l-1)\Vert_1}
+\hbox{O}\,((n_0+6l)^{-2})\, ,\cr
{e_A(2l)\over\Vert e_A(2l)\Vert_1} &= {e_A(2l-1)\over\Vert e_A(2l-1)\Vert_1}
+\hbox{O}\,((n_0+6l)^{-2})\, ,\cr
}
$$
and by applying the canonical involution one obtains similar formulas
at the order $2l+1$.
Therefore one can conclude that if there exist two vectors
$u_A$ and $u_D$ in $\R^\cA$ such that
$$
\eqalign{
\Vert u_A\Vert &=\Vert u_D\Vert =1\, , \cr
u_A &=\pmatrix{1/2 \cr 0\cr 0\cr 1/2\cr} +\hbox{O}\,(n_0^{-1})\, ,\cr
u_D &=\pmatrix{0 \cr 0\cr 1\cr 0\cr} +\hbox{O}\,(n_0^{-1})\, ,\cr
\lim_{l\rightarrow +\infty}{e_A(l)\over\Vert e_A(l)\Vert_1} &=
\lim_{l\rightarrow +\infty}{e_B(2l)\over\Vert e_B(2l)\Vert_1} =
\lim_{l\rightarrow +\infty}{e_C(2l)\over\Vert e_C(2l)\Vert_1} =
u_A\, , \cr
\lim_{l\rightarrow +\infty}{e_D(l)\over\Vert e_D(l)\Vert_1} &=
\lim_{l\rightarrow +\infty}{e_C(2l+1)\over\Vert e_C(2l+1)\Vert_1} =
\lim_{l\rightarrow +\infty}{e_B(2l+1)\over\Vert e_B(2l+1)\Vert_1} =
u_D\, .\cr
}
$$
It is now easy to see that each point $u$ of the segment $[u_A,u_D]
\subset(\R^+)^\cA$ is the lengths datum for an interval exchange map
with combinatorial datum $\pmatrix{ A & B & C & D\cr D & C & B & A\cr}$,
verifying Keane's condition and which is {\it not} uniquely ergodic:
the interval exchange maps of this one parameter family are
topologically conjugate.

 \vfill\eject
\noindent {\bf References}

\vskip .3 truecm \noindent
\item{[Ar]} P. Arnoux ``Ergodicit\'e g\'en\'erique des billiards
polygonaux [d'apr\`es Kerckhoff, Masur, Smillie]'' S\'eminaire
Bourbaki n. 696, Ast\'erisque 161--162, (1988), 203--221
\item{[Bo]} M. Boshernitzan ``A condition for minimal interval
exchange maps to be uniquely ergodic'' {\it Duke Math. J.} {\bf
52} (1985) 723--752
\item{[Ch]} Y. Cheung ``Hausdorff dimension of the set of nonergodic directions. With an appendix by M. Boshernitzan''
{\it Ann. of Math.} {\bf 158} (2003) 661--678.
\item{[Co]} J. Coffrey ``Some remarks concerning an example of a
minimal, non uniquely ergodic interval exchange map'' {\it Math.
Z.} {\bf 199} (1988) 577--580
\item{[FLP]} A. Fathi, F. Laudenbach and V. Poenaru
``Travaux de Thurston sur les surfaces'' Ast\'erisque 66-67
(1979).
\item{[Fo1]} G. Forni ``Solutions of the cohomological equation
for area-preserving flows on compact surfaces of higher genus''
{\it Annals of Mathematics} {\bf 146} (1997) 295-344.
\item{[Fo2]} G. Forni ``Deviation of ergodic averages for area-preserving
flows on surfaces of higher genus.'' {\it Annals of Mathematics}
{\bf 155 } (2002) 1--103.
\item{[Fo3]} G. Forni, private communication (2003)
\item{[GH]} W.H. Gottschalk, G.A. Hedlund ``Topological dynamics''
American Mathematical Society Colloquium Publications, {\bf 36.}
American Mathematical Society, Providence, R. I., (1955)
\item{[KH]} A. Katok and B. Hasselblatt ``Introduction to the modern
theory of dynamical systems'' Encyclopedia of Mathematics and its
Applications 54, Cambridge University Press, (1995).
\item{[KS]} A. Katok and A.M. Stepin ``Approximations in Ergodic
Theory'' {\it Russ. Math. Surv.} {\bf 22} (1967) 77--102
\item{[Ke1]} M. Keane ``Interval exchange transformations''
{\it Math. Z.} {\bf 141} (1975) 25--31
\item{[Ke2]} M. Keane ``Non--ergodic interval exchange transformations''
{\it Isr. J. Math.} {\bf 26} (1977) 188--196
\item{[Ker]} S. P. Kerckhoff ``Simplicial systems for interval
exchange maps and measured foliations'' {\it Ergod. Th. Dynam. Sys.}
{\bf 5} (1985) 257—271
\item{[KMS]} S. Kerckhoff, H. Masur and J. Smillie ``Ergodicity of
billiard flows and quadratic differentials'' {\it Ann. of Math.}
{\bf 124} (1986) 293--311
\item{[KN]} H. B. Keynes and D. Newton ``A ``Minimal'',
Non--Uniquely Ergodic Interval Exchange Transformation''
{\it Math. Z.} {\bf 148} (1976) 101--105
\item{[KR]}  M. Keane and G. Rauzy ``Stricte ergodicit\'e
des \'echanges d'intervalles'' {\it Math. Z.} {\bf 174} (1980)
203--212
\item{[KZ]} M. Kontsevich and A. Zorich
``Connected components of the moduli spaces of Abelian
differentials with prescribed singularities''
{\it Inv. Math.} {\bf 153} (2003) 631--678
\item{[Ma]} H. Masur ``Interval exchange transformations and measured
foliations'' {\it Annals of Mathematics} {\bf 115} (1982) 169--200
\item{[MMY]} S. Marmi, P. Moussa and J.--C. Yoccoz ``On the cohomological equation
for interval exchange maps'' {\it C. R. Math. Acad. Sci. Paris}
{\bf  336} (2003)  941--948
\item{[MS]} H. Masur and J. Smillie ``Quadratic differentials with prescribed
singularities and pseudo-Anosov diffeomorphisms'' {\it Comment. Math.
Helv.} {\bf 68} (1993) 289-307
\item{[Ra]} G. Rauzy ``\'Echanges d'intervalles et transformations
induites'' {\it Acta Arit.} (1979) 315--328
\item{[Re]} M. Rees ``An alternative approach to the ergodic
theory of measured foliations'' {\it Ergod. Th. Dyn. Sys.} {\bf 1}
(1981) 461--488
\item{[Ta]} S. Tabachnikov ``Billiards'' {\it Panoramas et
Synth\`eses, S.M.F.} {\bf 1} (1995)
\item{[V1]} W. Veech ``Interval exchange transformations''
{\it Journal d'Analyse Math\'ematique} {\bf 33} (1978) 222-272
\item{[V2]} W. Veech ``Gauss measures for transformations on the space of
interval exchange maps'' {\it Ann. of Math.} {\bf 115} (1982) 201--242
\item{[V3]} W. Veech ``The metric theory of interval exchange
transformations I. Generic spectral properties'' {\it Amer. J. of
Math.} {\bf 106} (1984) 1331--1359
\item{[V4]} W. Veech ``The metric theory of interval exchange
transformations II. Approximation by primitive interval exchanges'' {\it Amer. J. of
Math.} {\bf 106} (1984) 1361--1387
\item{[V5]} W. Veech ``The metric theory of interval exchange
transformations III. The Sah Arnoux Fathi invariant'' {\it Amer. J. of
Math.} {\bf 106} (1984) 1389--1421
\item{[V6]} W. Veech ``The Teichm\"uller geodesic flow'' Ann. of Math. 124 (1986)
441--530
\item{[V7]} W. Veech ``Moduli spaces of quadratic differentials'' Journal
d'Analyse Mathematique 55 (1990) 117--171
\item{[Y]} J.--C. Yoccoz ``Continued fraction algorithms for interval exchange maps: an
introduction'' preprint (2004), to appear in the proceedings of
conference Frontiers in Number Theory, Physics and Geometry, Les
Houches,  9 - 21 March 2003
\item{[Z1]} A. Zorich ``Finite Gauss measure on the space of interval
exchange transformations. Lyapunov exponents''
Annales de l'Institut Fourier Tome
46 fasc. 2 (1996) 325-370
\item{[Z2]} A. Zorich ``Deviation for interval exchange
transformations'' {\it Ergod. Th. Dyn. Sys.}{\bf 17} (1997),
1477--1499
\item{[Z3]} A. Zorich ``On Hyperplane Sections of Periodic Surfaces''
Amer. Math. Soc. Translations {\bf 179} (1997), 173--189
\item{[Z4]} A. Zorich ``How Do the Leaves of a Closed $1$--form Wind
Around a Surface?'' in Pseudoperiodic Topology, V. Arnold, M.
Kontsevich and A. Zorich editors, Amer. Math. Soc. Translations
{\bf 197} (1999) 135--178

\bye